\documentclass{article}
\usepackage{fullpage}
\usepackage{latexsym}
\usepackage{epsfig}
\usepackage{rotating}
\usepackage{amssymb}

\newcommand{\GL}{\mathbb{L}}
\newcommand{\GM}{\mathbb{M}}
\newcommand{\GK}{\mathbb{K}}

\newcommand{\GA}{\mathbb{A}}
\newcommand{\GS}{\mathbb{S}}

\newcommand{\GF}{\mathbb{F}}
\newcommand{\GB}{\mathbb{B}}

\newcommand{\GW}{\mathbb{W}}
\newcommand{\ns}{\mathbb{N}}
\newcommand{\zs}{\mathbb{Z}}
\newcommand{\qs}{\mathbb{Q}}
\newcommand{\rs}{\mathbb{R}}
\newcommand{\bx}{\bar x}
\newcommand{\by}{\bar y}

\newcommand{\cs}{\mathbb{C}}
\newcommand{\Z}{\mbox{\rm \lower0.3pt\hbox{$\angle\!\!\!$}Z}}

\newcommand{\Ref}[1]{(\ref{#1})}
\newcommand{\beq}{\begin{equation}}
\newcommand{\eeq}{\end{equation}}
\newcommand{\gf}{generating function}

\newcommand{\gfs}{generating functions}

\def\cqfd{\par\nopagebreak\rightline{\vrule height 3pt width 5pt depth 2pt}
\medbreak}

 \newtheorem{Theorem}{Theorem}
 \newtheorem{Propo}[Theorem]{Proposition}
\newtheorem{Coro}[Theorem]{Corollary}
\newtheorem{Lemma}[Theorem]{Lemma}

\begin{document}
\title{\bf
  Walks on the slit plane: other approaches}
\author{
 {\sc Mireille Bousquet-M{\'e}lou}\\ 
{\small CNRS, LaBRI, Universit\'e Bordeaux 1} \\
{\small 351 cours de la Lib\'eration} \\
{\small 33405 Talence Cedex,  France}\\
{\tt \small bousquet@labri.u-bordeaux.fr }
}
\date{}

\maketitle

\begin{abstract}
Let $\frak S$ be a finite subset of $\zs^2$. A walk on the slit plane
with steps in $\frak S$ is a sequence  $(0,0)=w_0, w_1, \ldots , w_n$  of
points of $\zs^2$ such that 
$w_{i+1}-w_i$ belongs to
$\frak S$ for all $i$, and none of the points $w_i$, $i\ge 1$, lie on
the half-line  ${\cal H}=\{(k,0): k\le 0\}$. 
In a recent paper, G.~Schaeffer and the author computed the length 
\gf \ $S(t)$ of walks on the slit plane for several sets $\frak S$.
All the \gfs \ thus obtained turned out to be algebraic: for instance,
on the ordinary square lattice, 
$$S(t) =\frac{(1+\sqrt{1+4t})^{1/2}(1+\sqrt{1-4t})^{1/2}}{2(1-4t)^{3/4}}. $$
The combinatorial reasons for this algebraicity remain obscure.

In this paper, we present two new approaches for solving slit plane
models. One of them simplifies and extends the functional equation
approach of the original paper. The other one is inspired by an argument
of Lawler; it is more combinatorial, and explains the algebraicity of
the product of three series related to the model. It can also be seen as
 an extension of the classical cycle lemma.
Both methods work for any set of steps $\frak S$.

We exhibit a large family of sets $\frak S$  for which the \gf \ of
walks on the slit plane
is algebraic, and
another family for which it 
is neither algebraic, nor even D-finite. These examples give a hint 
at where the border between algebraicity and transcendence lies, and
calls for a complete classification of the sets $\frak S$.
\end{abstract}

\section{Introduction}
\begin{figure}[ht]
\begin{center}
\epsfig{file=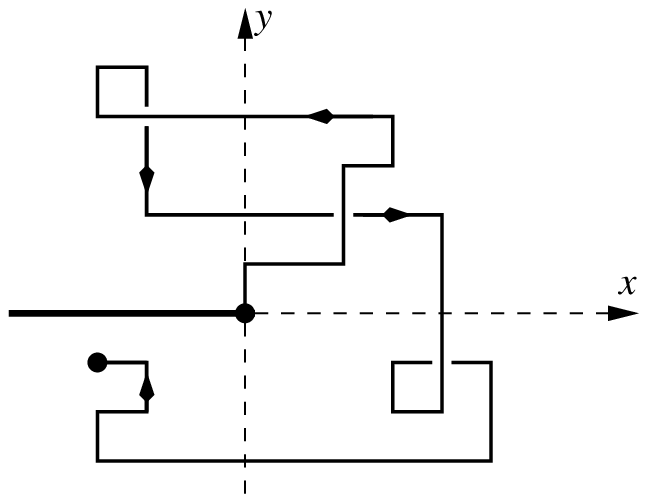}
\end{center}
\caption{A walk on the slit plane.}
\label{chemin}
\end{figure}

Let us consider square lattice walks that start from the origin $(0,0)$,
but never return to the 
horizontal half-axis ${\cal H}=\{(k,0): k\le 0\}$ once they have left their
starting point: we call them {\em walks on the slit plane\/}
(Figure~\ref{chemin}). We denote by $a_{i,j}(n)$  the number
 of such walks of length $n$ 
that end at $(i,j)$. In~\cite{prep-GS}, we proved that the associated
  \gf \ has a nice algebraic expression:
\begin{eqnarray}
S(x,y;t) &:=& \displaystyle 
\sum_{n\ge 0} t^n \sum_{(i,j) \in \zs ^2} 
a_{i,j}(n)x^iy^j \nonumber \\
&=&
\frac{\left(1-2t(1+\bar x)+\sqrt{1-4t}\right)^{1/2}
                   \left(1+2t(1-\bar x)+\sqrt{1+4t}\right)^{1/2}}
             {2(1-t(x+\bar x+y+\bar y))},\label{carre-sol}
 \end{eqnarray}
with $\bx =1/x$ and $\by=1/y$. We also studied the \gf \ for walks on the slit plane ending at a specific point $(i,j)$, proving, for instance, that
$a_{1,0}(2n+1) = C_{2n+1}$ and $a_{0,1}(2n+1)=4^n C_n$, where
$C_n={2n\choose n}/(n+1)$ is the $n$th Catalan number.
These two formulae had been conjectured by  Olivier Roques, 
 following a question raised by Rick Kenyon on a mailing list. 
The  enumeration of walks on the slit plane  was studied previously
in~\cite{considine,lawler},
mostly in probabilistic terms, and
is also related to~\cite{kesten}. However, all previously published 
results were asymptotic estimates, and to our knowledge,
it had never been realized that
the model could be exactly solved.

Algebraic \gfs \ and Catalan numbers call for ``classical''
--- if not bijective --- combinatorial proofs (see
\cite[Ex.~6.19]{stanley} for 66 other  
interpretations of these numbers). Indeed, the Florence group found a nice
bijective proof of the $C_{2n+1}$ result for walks on the slit plane
ending at $(1,0)$~\cite{firenze}. However, their approach so far only
works for this specific endpoint. As for the proof presented
in~\cite{prep-GS}, it is very atypical:
 from an obvious
combinatorial argument, we derived a  rather tricky functional equation,
which we solved.  This leaved us  with the
feeling that we had proved the identity~\Ref{carre-sol} without
understanding it 
completely. In particular, we have no combinatorial intuition as
why this series is algebraic. Our functional equation approach has,
nevertheless, also nice qualities: not only did it allow us to
prove~\Ref{carre-sol}, but it worked just as well for walks taking
their steps in a set $\frak S$ satisfying two simple conditions. 
For instance, we could solve the variation of the slit plane model
illustrated in Figure~\ref{turndiag}. 
For all these sets $\frak S$, we proved that the complete \gf \
$S(x,y;t)$ (the counterpart of~\Ref{carre-sol}) {\em and\/} the \gf \
for walks  ending at a prescribed position $(i,j)$ were always
algebraic.  This result has to be compared to the
case of ordinary walks (say, on the square lattice): their \gf \ is
rational, but, when restricted to walks ending at $(i,j)$, it becomes
transcendental. In other words, 
\begin{center}{\em forbidding a half-line adds some
algebraicity to the model\/}.
\end{center}
This rather mysterious algebraicity is the main motivation for
pursuing research on this topic: these results, we believe, need other
explanations.  

\begin{figure}[th]
\begin{center}
\epsfig{file=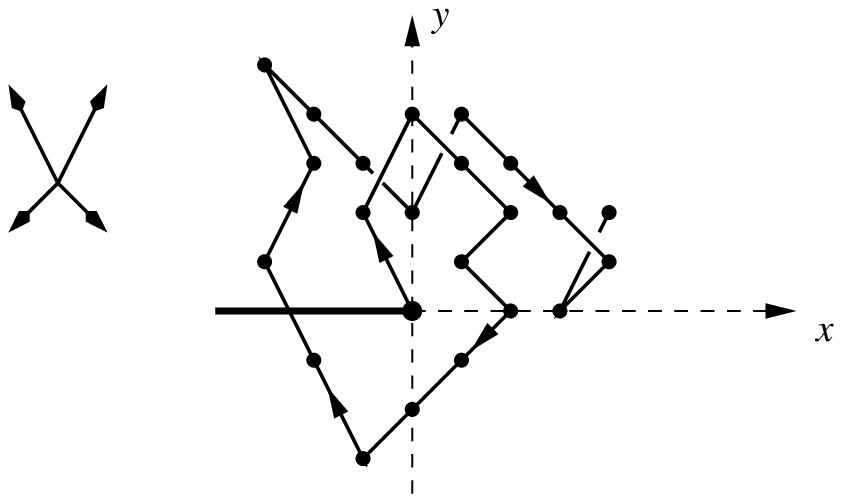}
\end{center}
\caption{A walk on the slit plane with steps in $\frak S=\{(1,2),(-1,2),(1,-1),(-1,-1)\}$.}
\label{chemin-trans}
\end{figure}

Let us confess right now that we shall not obtain here a completely
satisfactory answer. We shall, however, explain why, for any set of
steps $\frak S$, the {\em product\/} of three \gfs \ counting certain 
walks on the
slit plane is algebraic. We shall also clarify the limits of
algebraicity, proving that it holds for a larger class of steps than
in~\cite{prep-GS}, but that it fails for other simple sets $\frak S$, for
 which the \gf \ $S(x,y;t)$ is  not even 
D-finite\footnote{D-finite series, also called holonomic series,
 form a natural superset of algebraic series.}. 
These results rely on the fact that we are now able to solve slit
plane problems  {\em for any set $\frak S$\/}, in terms of a certain
factorization of the  \gf  \ for {\em bilateral walks\/} (walks ending
on the $x$-axis, but otherwise unconstrained). This \gf \ is
well-understood, and is always algebraic. Note that in these models,
we only forbid the walks to step 
on the half-line $\cal H$, not to cross it. Examples are given in
Figures~\ref{chemin-trans},~\ref{turndiag} and~\ref{chemin-triangle}.

\medskip
Let us now give a more detailed account of the contents  of this paper. In
Section~\ref{section-preliminaries}, we introduce an elementary terminology on
walks and words; we also recall some definitions and properties of
algebraic and D-finite power series. The power series fans can
certainly skip this part.

The paper really begins in Section~\ref{section-cycle}. We present a
``classical'' combinatorial proof of (some) slit plane results. 
This proof is based on the so-called  {\em cycle 
lemma\/}.
The underlying idea is
to recognize that the objects one wishes to count are in one-to-one
correspondence with cyclic conjugacy classes of some other objects,
which are easier to enumerate. 
The most classical application of the cycle lemma is the enumeration
of Dyck-like paths~\cite{motzkin,raney}.  
It is interesting to note that the florentine  proof of
the $C_{2n+1}$ formula for walks ending at $(1,0)$, although
completely different from our approach, is also based on conjugacy of
walks~\cite{firenze}. 
Our approach establishes a connection, valid for any set of steps
$\frak S$, between (some) walks on the slit plane and
bilateral walks.

However, the
cycle lemma suffers from a major drawback: it usually only works
for walks ending at a specific position. In our case, we can only
enumerate walks on the slit plane ending on the $x$-axis, as close as
possible to the origin. This means at the point $(1,0)$ for the models
of Figures~\ref{chemin} and~\ref{chemin-trans}, and at the point
$(2,0)$ for the diagonal model of Figure~\ref{turndiag}.

The next section, Section~\ref{section-lawler}, remedies this problem,
even though 
its starting point lies very far from the cycle lemma. It is inspired
by a ``very nonintuitive''
property that appears in Lawler's book, {\em Intersections of random
walks\/}~\cite{lawler}, to which Greg Lawler drew my attention.
This property, stated
in probabilistic terms, asserts that 
\begin{quote}
\begin{center}
{\em for a random walk starting
from $(0,0)$, with 
a killing rate\footnote{At each step, the walk is killed with
probability $\lambda$.} $\lambda$,  
the events \\``the walk avoids   $\{(k,0): k <0
\}$'' and 
``the walk avoids   $\{(k,0): k \ge 0
\}$'' \\are independent\/}. 
\end{center}
\end{quote}
From this result, Lawler derived that the
number of $n$-step walks on the slit plane, on the ordinary square
lattice, grows like $4^nn^{-1/4}$. We translate this property into
combinatorics, make a few elementary (but important) variations on it,
and end up with an identity that relates the \gf \ $S_0(x;t)$ for
 bilateral walks on the slit plane to the \gf \ $B(x;t)$ for
(ordinary) bilateral walks.
Again, this identity is valid for any set of steps $\frak S$. 
We prove that it contains and extends the cycle lemma result.

In Section~\ref{section-complete}, we prove that this identity
completely determines the \gf \ for walks on the slit plane: given any
set of steps, we are able to express $S_0(x;t)$ and $S(x,y;t)$ in
terms of the {\em 
canonical factorization\/} of the \gf \ $B(x;t)$ for bilateral walks. This kind
of factorization also played a major role in~\cite{prep-GS}, but for
more mysterious reasons; the identity ``\`a la Lawler'' makes it
completely natural. We argue that this expression of  $S(x,y;t)$
completely solves the problem: the series $\log S(x,y;t)$ is always
D-finite, and  a differential equation defining it can be obtained
in an algorithmic way.

In Section~\ref{section-functional}, we give an independent proof of
the connection between  $S(x,y;t)$  and the canonical factorization of
$B(x;t)$. We start from the very deep observation that an $n$-step
walk on the slit plane is obtained by adding a step to an $(n-1)$-step
walk on the slit plane. We translate this into a functional equation
for $S(x,y;t)$, and solve this equation (for any set of steps $\frak
S$) in a way that is so elementary that it looks fraudulous. This
approach is also the one that was used in the original
paper~\cite{prep-GS}: but the method we used for solving the
functional equation was 
more complicated, and only worked for certain sets $\frak S$.

The three previous sections provide tools for solving a specific
model, given a set of steps $\frak S$. They explain why $\log
S(x,y;t)$ is always D-finite. They  do not explain why, on the
square lattice, $S(x,y;t)$ is algebraic. The next two sections of the
paper are devoted to a partial 
exploration of the algebraic (and holonomic) properties of \gfs \ for
walks on the slit plane. In Section~\ref{section-algebraicity}, we
generalize the result of~\cite{prep-GS} on algebraicity 
to all sets of steps $\frak S$
having {\em small height variations\/}: by this, we mean that adding
a step to a walk  never modifies its ordinate by more than one unit.
This condition prevents a
walk from crossing the forbidden half-line without hitting it. We
solve explicitly three examples: the ordinary square lattice of
Figure~\ref{chemin}, the diagonal square lattice of
Figure~\ref{turndiag} (both already solved in~\cite{prep-GS}), and the
triangular lattice of Figure~\ref{chemin-triangle}. Would
$S(x,y;t)$ be algebraic for {\em all\/} sets of steps? No: we study in
Section~\ref{section-transcendental} some extensions of the diagonal
square lattice  for which $S(x,y;t)$ is neither algebraic, nor
even D-finite. An  example is given in
Figure~\ref{chemin-trans}.

We conclude the paper by mentioning possible directions of further
research. This includes a complete
classification of the sets $\frak S$ according to the
algebraicity/transcendence of the corresponding \gf : an obvious
conjecture is that algebraicity holds if and only if a walk cannot
cross the forbidden half-line without hitting it. This also includes
higher-dimensional models, to which the approach of this paper can be
applied.
\section{Preliminaries}
\label{section-preliminaries}

\subsection{Walks on the slit plane: definitions and notations} 
\label{definitions}
Let $\frak S$ be a finite subset of $\zs ^2$. A walk with steps in
$\frak S$ is a
finite sequence $w = (w_0, w_1, \ldots , w_n)$  of vertices of $\zs
^2$ such that $w_0=(0,0)$ and $w_{i}-w_{i-1} \in \frak S$ for $1 \le i \le
n$. The number of 
steps, $n$, is the {\em length\/} of $w$. The {\em endpoint\/}  of $w$ is
$w_n$, and its abscissa is denoted $\delta(w)$.  The \gf \ for  a set $\cal A$ of walks is the formal series
$$ A(x,y;t)=\sum_{n \ge 0}t^n \sum_{i,j \in \zs} a_{i,j}(n) x^i 
y^j ,$$ where $a_{i,j}(n)$ is the number  of walks of  $\cal A$
that have length $n$  and  end at $(i,j)$. This series is a formal
power series in $t$ whose coefficients are polynomials in $x, y, 1/x,
1/y$. We shall often denote $\bx =1/x$ and $\by =1/y$.

We say that the walk  $w$ avoids the 
half-line ${\cal H}=\{(k,0), k\le 0\}$ if none of the vertices $w_1, \ldots ,
w_n$ belong to $\cal H$. We call $w$ a {\em walk on the slit plane\/}
(Figure~\ref{chemin-trans}). We denote by $\cal S$ the set of walks on the slit
plane, and by $S(x,y;t)$ the corresponding \gf . For $j \in \zs$, we
denote by  ${\cal S}_j$ the set of walks on the slit
plane ending at ordinate $j$, and by $S_j(x;t)$ the corresponding
(two-variable) \gf . Finally, for  $i, j \in \zs$, we
denote by  ${\cal S}_{i,j}$ the set of walks on the slit
plane ending at the point $(i,j)$, and by $S_{i,j}(t)$ the corresponding
(one-variable) \gf .

We say that the walk  $w= (w_0, w_1, \ldots , w_n)$ is a {\em loop\/} if
$w_0=w_n=(0,0)$  and  none of the vertices $w_i$ belong to the negative
$x$-axis $\{(k,0), k< 0\}$. 

A walk with steps in $\frak S$ is called {\em bilateral\/} if it ends
on the $x$-axis.

\subsection{Walks and words} 
By definition, a walk is a sequence of steps of $\frak S$. Hence, it
will be convenient to consider walks as words on the alphabet $\frak
S$. The set of words $\frak S^*$ is equipped with the usual
concatenation product. The empty word is denoted $\epsilon$.
The number of occurrences of the letter $a$ in the word $w$ is
denoted $|w|_a$. Any word $w$ of   $\frak S^*$  will be thought of as
a walk  starting 
from $(0,0)$, and it will be convenient to make no distinction 
between the walk and the word, using expressions like ``the endpoint
of $w$'', ``the number of east steps in $w$'', etc.

A {\em language\/} is a subset   ${\cal A}$ of $ \frak S^*$; in other
words, a set of walks. 
We denote by $\GA$ the {\em non-commutative \gf \/} of
$\cal A$, that is, the formal power series
$$\GA = \sum_{w \in {\cal A}} w.$$
The product ${\cal A}_1{\cal A}_2$ of two languages ${\cal A}_1$ and
${\cal A}_2$ is the set 
of words $w$ that can be written as $uv$, where $u \in {\cal A}_1$ and
$v \in {\cal A}_2$. The product is {\em non-ambiguous\/} if any word in the
product language can be factored in a unique way as  $uv$, with $u \in
{\cal A}_1$ and $v \in {\cal A}_2$. In this case, the non-commutative
\gf \ for the product  ${\cal A}_1{\cal A}_2$ is  $\GA_1 \GA_2$.


\subsection{Classes of formal power series} \label{rappels-DF}

 Given a ring $\GL$ and $k$ indeterminates
$x_1, \ldots , x_k$, we denote by
 $\GL[x_1, \ldots , x_k]$ the ring of polynomials in $x_1,
\ldots , x_k$ with coefficients in $\GL$,
and by $\GL[[x_1, \ldots , x_k]]$ the ring of {\em formal power
series\/} in $x_1, 
\ldots , x_k$ with coefficients in $\GL$, that is, formal sums
\beq \sum_{n_1\ge 0, \ldots , n_k \ge 0} a_{n_1, \ldots ,
n_k}x_1^{n_1}\cdots x_k^{n_k}, \label{formal-series}\eeq  
where $ a_{n_1, \ldots , n_k} \in \GL$.
 If $\GL$ is a field, we denote by
 $\GL(x_1, \ldots , x_k)$ the field of rational functions in $x_1,
\ldots , x_k$ with coefficients in $\GL$.
A {\em Laurent polynomial\/} in the $x_i$ is a polynomial in the
$x_i$ and the $\bx _i=1/x_i$. A {\em Laurent series\/} in the $x_i$ is
a series of the form~\Ref{formal-series} in which the summation runs
over ${n_1\ge N_1, \ldots , n_k \ge N_k}$ for $N_1, \ldots , N_k \in \zs$.

Let $F$ be a series of $\GL[[x_1, \ldots , x_k]]$. Let $i \le k$  and write
$$F= \sum_{ n \ge 0} F_{n}\ x_i^n$$
where the coefficients $F_{n}$ are power series in the remaining
variables $x_\ell$, $\ell \not = i$. We denote by $[x_i^n]F$ the coefficient
of $x_i^n$ in $F$, that is,
$[x_i^n]F=F_n.$
 Let $i,j \le k$  and write
$$F= \sum_{m, n \ge 0} F_{m,n}\ x_i^m x_j ^n$$
where the coefficients $F_{m,n}$ are power series in the 
variables $x_\ell$, $\ell \not = i,j$. The {\em diagonal\/}  of $F$ with respect to
$x_i$ and $x_j$ is defined to be  $\sum _{n \ge 0}F_{n,n} x_i^n $.  

\medskip
Assume $\GL$ is a field. A series $F$ in $\GL[[x_1, \ldots , x_k]]$ is
{\em rational\/} if 
there exist non-trivial polynomials $P$ and $Q$ in $\GL[x_1, \ldots ,
x_k]$ such that $PF=Q$. It is {\em algebraic\/}
(over the field  $\GL(x_1, \ldots , x_k)$) if
it there exists a non-trivial polynomial $P$ with coefficients in
$\GL$ such that 
$P(F,x_1, \ldots , x_k)=0.$ The sum and product of algebraic series
is algebraic.  The diagonals of a rational series  are
algebraic~\cite[p.~179]{stanley}.

The series $F$  is {\em D-finite\/} (or {\em holonomic\/})  
if the partial derivatives of $F$ span a finite 
dimensional vector space on the field $\GL(x_1, \ldots , x_k)$
(see~\cite{stanleyDF} for the one-variable case,
and~\cite{lipshitz-diag,lipshitz-df} otherwise).  In other 
words, for $1\le 
i\le k$, the series $F$ satisfies a non-trivial partial differential
equation of the form
$$\sum_{\ell=0}^{d_i}P_{\ell,i}\frac{\partial ^\ell
F}{\partial x_i^\ell} =0,$$
where $P_{\ell,i}$ is a polynomial in the $x_j$.
Any algebraic series is holonomic. The
sum and product of two 
holonomic series is still holonomic. The
specializations of an holonomic series (obtained by giving 
values from $\GL$  to some of the variables) are holonomic, if
well-defined. Moreover, if $F$ is an {\em algebraic\/} series and 
$G(t)$ is a holonomic series of one variable, then the substitution
$G(F)$ (if well-defined) is
holonomic~\cite[Prop.~2.3]{lipshitz-df}. Most importantly, if $F$ is 
holonomic, then any diagonal of $F$  is also
D-finite~\cite{lipshitz-diag}.

\section{Counting (some) walks on the slit plane by the cycle lemma}
\label{section-cycle}
  
\subsection{The ordinary square lattice} 
We focus here on the ordinary square lattice: the set of steps is
$\frak S=\{(0,1), (1,0), (0,-1), (-1,0)\}=\{o,n,e,s\}$,  the letter $o$ 
(resp. $n,e,s$) standing for a west (resp. north, east, south)
step\footnote{The french word for {\em west\/} being {\em
ouest}.}. Recall that a {\em bilateral walk\/}
is a walk that ends on the $x$-axis. A (bicolored) {\em Motzkin
walk\/} is a bilateral walk that never visits a point with a negative
ordinate. We shall enumerate these walks not only by their length
 and endpoint (variables $t$ and $x$), but also by   the number of
their vertical 
steps, using an additional indeterminate $v$.
The associated \gfs \  can be determined
by a standard argument. 
\begin{Lemma}\label{lemmeMB}
The \gf \ for bicolored Motzkin walks is
\begin{eqnarray*}
M (x;t,v) 
&=& \frac{1-t(x+\bar x) -\sqrt{[1-t(x+\bar x
+2v)][1-t(x+\bar x -2v)]} }{2t^2v^2}.
\end{eqnarray*}
The \gf \ for bilateral walks is
$$
B (x;t,v) =  \frac{1}{\sqrt{[1-t(x+\bar x
+2v)][1-t(x+\bar x -2v)]}}.
$$
\end{Lemma}
{\bf Proof.} 
The simplest approach consists in factoring (non-empty) Motzkin walks
at the first place where they return to the horizontal axis. This
gives, for their
non-commutative generating function,
$$\GM = \epsilon +o\  \GM + e\ \GM + n\  \GM \ s\  \GM,$$
so that
$$M(x;t,v)=1+ t(x+\bar x) M(x;t,v) + t^2v^2M(x;t,v)^2.$$ The
expression of $M(x;t,v)$ follows. The same principle, applied to  bilateral
walks, gives
$$\GB=  \epsilon +o\  \GB + e\ \GB + n\  \GM \ s\  \GB + s \ \Phi(\GM)\  n
\ \GB,$$ 
where $\Phi$ is the morphism that  flips the walks around the
horizontal axis: $\Phi(n)=s, \Phi(s)=n,
\Phi(e)=e$ and $\Phi(o)=o$. Hence
$$B(x;t,v)=1+t(x+\bx) B(x;t,v) + 2t^2v^2  B(x;t,v) M(x;t,v).$$
The expression of $B(x;t,v)$ follows. 
\cqfd

Recall that the $n$th Catalan number $C_n$ is defined by:
$$C_n=\frac 1 {n+1} {{2n} \choose n}.$$
\begin{Propo} \label{propo-cycle-carre}
 The number of walks of length $2n+1$ on the slit plane  ending at
 $(1,0)$ is  the  Catalan number $ C_{2n+1}$.
More precisely,  the number of such walks having exactly $2m$ 
vertical steps 
is $$
4^{m}{{2n } \choose {2m}}C_{n-m}.$$  
This is also the number of bicolored Motzkin  walks of length $2n$ with
$2m$  horizontal steps.
\end{Propo}
{\bf Proof.}
Let us say that a non-empty bilateral walk $w$ is {\em primitive\/} if the
only vertices of $w$ that lie 
on the $x$-axis are its endpoints. We denote by ${\cal P}$ the set of
primitive words, and by $P(x;t,v)$ their \gf .

\begin{figure}[ht]
\begin{center}
\input{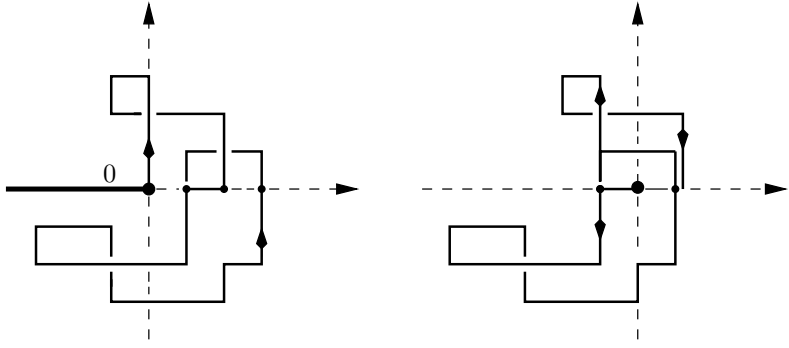}
\end{center}
\caption{A word  $w=w_1 \cdots w_{4}$ of ${\cal S}_{1,0}$ 
and one of its conjugates $v=w_2 w_3w_4w_1=v_1v_2v_3v_4$.}
\label{cyclic-horizontal}
\end{figure}

   Let $w \in {\cal S}_{1,0}$, and
assume that exactly $(k+1)$ of its vertices lie on the $x$-axis. Then
$w$ can be written in a unique way as the concatenation $w_1w_2 \cdots
w_k$ of $k$ primitive words. For $i \in [k]$, let $w^{(i)}$ be
obtained by permuting cyclically the factors $w_\ell$, starting from
$w_i$~(see Figure~\ref{cyclic-horizontal}). That is, $w^{(i)}= w_i
w_{i+1} \cdots w_kw_1 \cdots 
w_{i-1}$. Note that $\delta (w^{(i)}) =\delta (w)=1$.

Conversely, let $v = v_1 \cdots v_k$ be a word of ${\cal P}^k$
such that $ \delta(v) =1$. Let $v_j$ be the primitive factor following  the 
last vertex of $v$ that lies on the $x$-axis and has a minimal
abscissa. Clearly,  $v^{(j)}:=v_j
v_{j+1} \cdots v_{k} v_1 \cdots v_{j-1}$ belongs to $ {\cal
S}_{1,0}$. We claim that no other word $v^{(\ell)}$ belongs to this
set. Indeed, if $\ell <j$, then  the prefix $v_\ell \cdots v_{j -1}$ 
of  $v^{(\ell)}$  ends at abscissa $\delta(v_1\cdots
v_{j -1})- 
\delta(v_1\cdots v_{\ell -1})$. By definition of
$j$, this abscissa  is non-positive, so that this prefix of
$v^{(\ell)}$ ends on $\cal H$. Similarly, if $\ell>j$, then $v_\ell \cdots v_k v_1 \cdots
v_{j -1}$ is a prefix of  $v^{(\ell)}$ that ends  at abscissa $1-
\delta(v_1\cdots v_{\ell -1}) +\delta(v_1\cdots v_{j -1}) <1$, that
is, on the forbidden half-line $\cal H$.

This proves that the map $(w,i) \mapsto w^{(i)}$ establishes a
 one-to-one correspondence between:


-- pairs $(w,i)$, where $w$ is a walk  of ${\cal S}_{1,0}$ formed of $k$ primitive factors and $i\in[1,k]$, and

-- bilateral walks $v$  ending at
 $(1,0)$ and formed of $k$ primitive factors.

\noindent This bijection preserves the number of vertical/horizontal
 steps. Consequently, the \gf \ for walks  of
${\cal S}_{1,0}$ formed of $k$ primitive factors
is $$ \frac{1}{k}[x] P(x;t,v)^k,$$ and summing over $k \ge 1$ gives
$$S_{1,0}(t,v) =  [x] \sum_{k \ge 1} \frac{ P(x;t,v)^k}{k} =  [x] \log
\frac{1}{1-P(x;t,v)} =  [x] \log B(x;t,v),$$
because a bilateral walk is simply a sequence of primitive walks. The
explicit value of $B(x;t,v)$, given in Lemma~\ref{lemmeMB}, yields
\begin{eqnarray}
S_{1,0}(t,v)
&=& \frac{1}{2} \ [x] \left[ 
\log \frac{1}{1-t(x+\bar x +2v)} + 
\log \frac{1}{1-t(x+\bar x -2v)} \right] \nonumber \\
&=& \frac{1}{2} \ [x] \sum_{n \ge 1} \frac{t^n}{n} 
\left[ (x+ \bar x +2v)^n + (x+\bar x -2v)^n\right]\nonumber \\
&=& [x] \sum_{n \ge 1} \frac{t^n}{n} 
\sum_{m \ge 0} (2v)^{2m} {n \choose 2m}  (x + \bar x)^{n-2m}\nonumber \\
&=&\sum_{n \ge 0} \frac{t^{2n+1}}{2n+1} \sum_{m=0}^n (2v)^{2m} {{2n+1}
\choose {2m}} {{2n-2m+1} \choose {n-m}} .\nonumber 
\end{eqnarray}
This gives the second result of Proposition~\ref{propo-cycle-carre}.
Now when $v=1$,
$$B(x;t,1)= \frac{1}{\sqrt{(1-t\bx(1+x)^2)(1-t\bx(1-x)^2)}}$$
and a similar calculation yields the first result of the proposition.
\cqfd

\medskip
\noindent{\bf Remark.} It is easy to see that $4^m {{2n} \choose {2m}}
C_{n-m}$ is the number of bicolored Motzkin words of length $2n$ with
$2m$ horizontal steps: such words are obtained by shuffling a word of
length $2m$ on the alphabet $\{e,o\}$ with a Dyck word of length
$2n-2m$ on the alphabet  $\{n,s\}$. And Dyck words are well-known to
be counted by Catalan numbers. A nice bijection between walks of
${\cal S}_{1,0}$ and  bicolored Motzkin words proving
Proposition~\ref{propo-cycle-carre} is presented in~\cite{firenze}.

\subsection{Walks with steps in an arbitrary set $\frak S$ } 
We first focus on bilateral walks:   whatever
the set $\frak S$, their \gf \ is  algebraic.
\begin{Lemma} \label{B-general}
Let $K(x,y;t)$ be the following polynomial:
\beq K(x,y;t)=1-t\sum_{(i,j) \in \frak S} x^i y^j.\label{kernel}\eeq
The \gf \  for bilateral walks is algebraic, given by
$$ B(x;t) = [y^0] \frac{1}{K( x,y;t)}. $$
\end{Lemma}
{\bf Proof.}   A walk being a sequence of steps, the \gf \ $W(x,y;t)$ for all 
 walks with steps in $\frak S$ is the  rational function
$W(x,y;t)= {1}/{K(x,y;t)}.$ The expression of $B$ is a direct
 consequence of the definition of bilateral walks. 

 Let us now discuss the algebraic nature of $B$. The idea is that
taking a constant term is equivalent to taking a diagonal\footnote{An
argument that goes against nature: several proofs that the diagonal of
a rational series is algebraic start by transforming the diagonal
problem into a constant term one...}. We
introduce two new variables, $\bx$ and $\by$, independent, for once,
of $x$ and $y$. We refine the weight of the walks by giving to a
step $(i,j)$ the weight $tx^iy^j$ if $i\ge 0$ and $j\ge 0$,  the weight
$t\bx^{-i}y^j$ if $i\le 0$ and $j\ge 0$, etc. Let 
$\overline W(x,\bx,y,\by;t)$ be the \gf \  for all walks, in which each
walk is weighted by the product of the weights of its steps. This
series  is still rational, and
$B(x;t)$ is derived from $\overline W$ first by taking the diagonal with
respect to  $y$
and $\by$ (which yields an algebraic series), then by setting $y=1$
and $\bx =1/x$ (which preserves algebraicity).
\cqfd

We now state the generalization of Proposition~\ref{propo-cycle-carre}.
\begin{Propo}[The cycle lemma result]\label{propo-cycle-1}
Let $p$ be the smallest positive integer such that there exists a walk
on the slit plane ending at the point $(p,0)$. 
Then the \gf \ for such walks is
$$ S_{p,0}(t)=  [x^p] \log B(x; t) $$
where $B(x;t)$ counts bilateral walks. This series is D-finite.
\end{Propo}
{\bf Proof.} The proof mimics that of
Proposition~\ref{propo-cycle-carre}, with one subtlety.
Again, we say that a non-empty bilateral walk $w$ is primitive if 
the only vertices of $w$ that lie 
on the $x$-axis are its endpoints. 

Let $v = v_1 \cdots v_k$ be a bilateral word formed of $k$ primitive
factors, and 
such that $ \delta(v) =m>0$. Let $v_j$ be the primitive factor following  the 
last vertex of $v$ that lies on the $x$-axis and has a minimal
abscissa. Clearly,  $v^{(j)}=v_j
v_{j+1} \cdots v_{k} v_1 \cdots v_{j-1}$ belongs to $ {\cal
S}_{m,0}$.  In particular, $m \ge p$. This shows that {\em a
bilateral walk 
ending at a positive abscissa actually ends to the right of $(p,0)$}.

We can now restrict the study to the case where the endpoint of $v$ is
exactly $(p,0)$. We claim that no other word $v^{(\ell)}$ belongs to  $ {\cal
S}_{p,0}$. Indeed, if $\ell <j$, then  the prefix $v_\ell \cdots v_{j -1}$ 
of  $v^{(\ell)}$  ends at abscissa $\delta(v_1\cdots
v_{j -1})- \delta(v_1\cdots v_{\ell -1})$. By definition of
$j$, this abscissa  is non-positive, so that this prefix of
$v^{(\ell)}$ ends on $\cal H$. Similarly, if $\ell>j$, then $v_\ell
\cdots v_k v_1 \cdots 
v_{j -1}$ is a prefix of  $v^{(\ell)}$ that ends  at abscissa $p-
\delta(v_1\cdots v_{\ell -1}) +\delta(v_1\cdots v_{j -1}) <p$. But we
have just proved that this forces this prefix (a bilateral walk) to
end at a non-positive abscissa, hence on $\cal H$. 

From this point, we conclude exactly as in the proof of
Proposition~\ref{propo-cycle-carre} to obtain the expression of $S_{p,0}$.

\medskip
Let us now discuss the holonomic nature of this series. As in the proof of
Lemma~\ref{B-general}, we refine the weight of the steps by
introducing two new variables $\bx$ and $\by$, independent of $x$
and $y$. The \gf \  $\overline B(x, \bx;t)$ that counts bilateral walks is
the diagonal of a rational series, and thus is an algebraic series in
$x$, $\bx$ and $t$. Hence its logarithm is 
D-finite. Finally,  $S_{p,0}(t)$ is obtained by first taking the
diagonal of the series $\bx ^p \log B(x, \bx ;t)$ with respect to  $x$
and $\bx$, which yields a D-finite series, and then by setting $x=1$,
which preserves holonomy.
\cqfd 

\medskip
The fact that we are dealing with two-dimensional walks should not
hide that our argument  is essentially one-dimensional. Let us state
explicitly the result for walks  on $\zs$ that underlies
Proposition~\ref{propo-cycle-1}. 
\begin{Lemma}\label{1D-cycle}
Let  $\frak P$ be a set of steps in $\zs$, and let ${\bf
\lambda}=(\lambda_i)_{i\in\frak P}$ be a sequence of variables describing
the weights of these steps.
We define the weight  
of a path $w$ with steps in $\frak P$ to be
the product of the weights of its steps. Let $p$ be the smallest
positive integer such that there exists a walk ending at $p$. Then the
\gf \ for $n$-step walks that start from $0$, end at $p$ and always
have a positive level once they have left their starting point is
$$\frac 1 n [x^p] \left( \sum_{i \in \frak P} \lambda_i
x^i\right)^n.$$
Hence the \gf \ for walks ending at $p$ is
$$S_p=[x^p] \log \frac 1 {1-\sum_{i \in \frak P} \lambda_i
x^i}.$$
\end{Lemma}

\subsection{The diagonal square lattice} 
\begin{figure}[ht]
\begin{center}
\epsfig{file=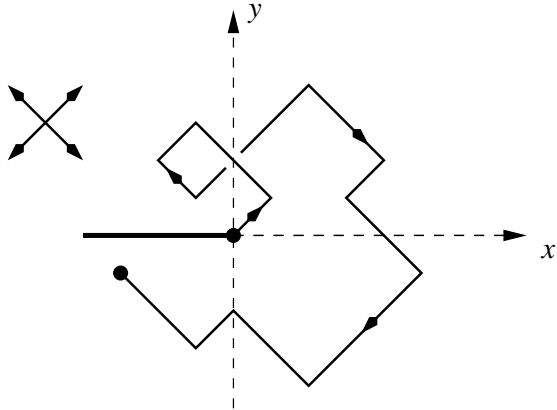}
\end{center}
\caption{A walk  on the slit plane with  diagonal steps.}
\label{turndiag}
\end{figure}

We now apply the general result of Proposition~\ref{propo-cycle-1} to
the diagonal 
version of the square lattice, illustrated in Figure~\ref{turndiag}.   
The set of steps is $\frak S=\{(1,1),(1,-1),(-1,1),(-1,-1)\}$.
By analogy with what we did on the square lattice, we
 put an additional weight $v$ on the ``vertical'' steps $(1,1)$
and $(-1,-1)$. This, we shall see,  gives a new occurrence of the
Narayana  distribution~\cite[p.~273]{narayana,stanley}. 

\begin{Lemma}\label{lemmaMBdiag}
The \gf \ for bilateral walks on the diagonal square lattice is
$$B(x;t,v)=\frac{1}{\sqrt{1-4t^2(v+x^2)(v+\bar x^2)}}.$$
\end{Lemma}
{\bf Proof.} The principle is the same as in
Lemma~\ref{lemmeMB}. The \gf \ $M(x;t,v)$ for bilateral walks that
stay on or above the $x$-axis satisfies
$$M(x;t,v)=1+t(xv+\bx)M(x;t,v)t(x+\bx
v)M(x;t,v)=1+t^2(v+x^2)(v+\bx^2)M(x;t,v)^2.$$ 
 The  \gf \ $B$ for bilateral walks is related to $M$ by
$$B(x;t,v)=1+2t^2(v+x^2)(v+\bx^2)M(x;t,v)B(x;t,v).$$
(The factor $2$ comes from the fact that primitive bilateral words
starting with a down step are obtained by reading from right to left
primitive bilateral words starting with an up step.) The results follows.
\cqfd
\begin{Propo}\label{propo-cycle-diagonal}
 The number of walks of length $2n$  with diagonal steps going from
 $(0,0)$ to $(2,0)$ on  the  slit plane is
$$
\frac{1}{2} \ 4^n C_{n}
$$
where $C_n$ denotes the $n$th Catalan number.
More precisely,  the number of such walks having exactly $2m-1$ 
 steps in the set $\{(-1,-1), (1,1)\}$
is 
$$
\frac{1}{2} \ \frac{4^n}{n} {n \choose m} {n \choose {m-1}}.$$
The underlying  refinement of Catalan numbers is known as the Narayana
distribution. 
\end{Propo}
{\bf Proof.} The proof is an immediate combination of 
Proposition~\ref{propo-cycle-1} (with $p=2$) and Lemma~\ref{lemmaMBdiag}:
\begin{eqnarray*}
  S_{2,0}(t,v) 
&=&
[x^2] \log {B(x;t,v)}  \\
&=&[x^2] \log \frac{1}{\sqrt{1-4t^2(v+x^2)(v+\bar x^2)}}  \\
&=& \frac{1}{2} [x^2]  \sum_{n \ge 1} \frac{4^n t^{2n}}{n}
(v+x^2) ^n(v+\bar x^2)^n\\
&=& \frac{1}{2} \sum_{n \ge 1} \frac{4^n t^{2n}}{n} \sum_{m=1}^n
v^{2m-1} {n \choose m} {n \choose {m-1}}.
\end{eqnarray*}
This gives the second result of Proposition~\ref{propo-cycle-diagonal}.
Now when $v=1$,
$$B(x;t,1)= \frac{1}{\sqrt{1-4t^2\bx^2(1+x^2)^2}}$$
and a similar calculation yields the first result of the proposition.
\cqfd
\section{On the steps of Greg Lawler}
\label{section-lawler}

We consider walks on the slit plane with steps in a generic set $\frak
S$. We denote by $\overline \frak S$ the set of steps obtained by
reversing the steps of $\frak S$. That is, $\overline \frak S = \{
(-i,-j) : (i,j) \in \frak S\}$. Accordingly, the \gfs \ for walks on
the slit plane with
steps in $\overline \frak S$ are denoted $\bar S(x,y;t)$, $\bar
S_0(x;t)$, and so on. Observe that the \gfs \ for loops with steps in  
$\frak S$ and  $\overline \frak S$ are equal.

\begin{figure}[ht]
\begin{center}
\input{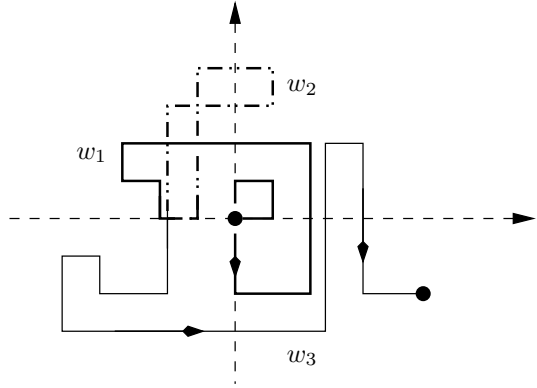}
\end{center}
\caption{Factorisation of an ordinary walk.}
\label{gilles}
\end{figure}

The following proposition relates a number of ``difficult'' series
(counting walks on the slit plane) to the \gf \ $B(x;t)$ of bilateral
walks, which is simpler to determine (see Lemma~\ref{B-general}). We
shall see that it contains and extends the cycle lemma result. Note that
a {\em bilateral\/} walk on the slit plane is simply a walk of ${\cal
S}_0$.
\begin{Propo}\label{greg}
The \gfs \ for bilateral  walks on the slit plane with steps in $\frak
S$ (resp.~$\overline \frak S$) and the \gf \ $L(t)$ for loops are
related by
\beq \bar{S}_0(\bx;t) L(t) S_0(x;t)=B(x;t),\label{greg-equation}\eeq
where $B(x;t)$ counts bilateral walks. The \gf \ for walks on the slit
plane (with steps in $\frak S$) can be expressed in terms of the \gf \
for {\em bilateral}  walks on the slit plane:
$$S(x,y;t)= \frac{S_0(x;t)}{K(x,y;t)B(x;t)},$$
where $K(x,y;t)$ is the polynomial given by~{\em\Ref{kernel}}. 
\end{Propo}
{\bf Proof.}
Let $w$ be an ordinary walk with steps in $\frak S$. Let $(i,0)$ be
the leftmost point of 
the half-line $\cal H$ visited by $w$. By looking at the first and
last visit of $w$ 
to this point, we  factor $w$ into three
walks $w_1, w_2, w_3$, where $w_1$ is obtained by reversing a bilateral walk on
the slit plane with steps in  $\overline \frak S$, $w_2$ is a loop and
$w_3$ belongs to  ${\cal  S}$ (Figure~\ref{gilles}). This gives 
\beq \frac 1 {K(x,y;t)} = \bar{S}_0(\bx;t) L(t) S(x,y;t).\label{greg-ordinary}\eeq
By restricting the argument to {\em bilateral\/} walks, we obtain the
first identity of the proposition.  Finally, eliminating
$\bar{S}_0(\bx;t)$ between~\Ref{greg-equation} and~\Ref{greg-ordinary}
gives the  second identity of Proposition~\ref{greg}.

\cqfd

Let us now explain why we claim that the first part of
Proposition~\ref{greg} generalizes the cycle lemma result of
Proposition~\ref{propo-cycle-1}. Let us take logarithms in~\Ref{greg-equation}:
\beq \log \bar{S}_0(\bx;t) + \log L(t) + \log S_0(x;t)=\log
B(x;t). \label{log-lawler}\eeq
The only walk of ${\cal S}_0$ ending at abscissa $0$ is the empty
walk. Hence 
$$S_0(x;t)= 1+\sum_{i\ge 1} S_{i,0}(t) x^i .$$ For $i\ge 1$, let us
extract the coefficient of $x^i$ from~\Ref{log-lawler}. The series 
$\log \bar{S}_0(\bx;t)$ and $ \log L(t) $ do not contribute, and we
obtain the following proposition.
\begin{Propo}\label{propo-cycle-2}
For $i \ge 1$, the  \gf \ $S_{i,0}(t)$ for walks on the slit plane ending
at $(i,0)$ is D-finite and can be computed by induction on $i$ via
the following identity: 
$$\sum_{k =1}^i\frac{(-1)^{k-1}}{k} \sum_{i_1+\cdots+i_k=i\atop i_1
>0, \ldots , i_k >0}
S_{i_1,0}\cdots S_{i_k,0} =  [x^{i}] \log B(x;t)$$
where $B(x;t)$ counts bilateral walks.
\end{Propo}
Let $p$ be the smallest positive integer such that there exists a walk
on the slit plane ending at the point $(p,0)$. The case $i=p$ of the
above proposition is exactly Proposition~\ref{propo-cycle-1}. The
argument that proves the holonomy of  $S_{i,0}(t)$ for general $i$ is
the same as in Proposition~\ref{propo-cycle-1}. 

\medskip Let us illustrate the above result by the example of the
ordinary square lattice. For this model, $p=1$, and the first three
instances of Proposition~\ref{propo-cycle-2} are:
\begin{eqnarray*}
S_{1,0} &=& [x] \log B(x;t),\\
\displaystyle S_{2,0}-\frac {1}{2} S_{1,0}^2 &=& [x^{2}] \log B(x;t),\\
\displaystyle S_{3,0}- S_{1,0}S_{2,0} + \frac {1}{3} S_{1,0}^3 &=&
[x^{3}] \log B(x;t).
\end{eqnarray*}
The  \gf \ $B(x;t)$ for bilateral walks on the square
lattice is derived from Lemma~\ref{lemmeMB},
$$B(x;t)= \frac{1}{\sqrt{(1-t\bx(1+x)^2)(1-t\bx(1-x)^2)}},$$
and allows us to compute explicitly
$$[x^{i}] \log B(x;t)= \sum_{n\ge 0}\frac{t^{2n+i}}{2n+i}
{{4n+2i}\choose{2n}}.$$
This series can actually be seen to be algebraic, so that $S_{i,0}(t)$
is not simply D-finite, but also algebraic. We shall exhibit in
Section~\ref{section-algebraicity} a whole family of sets $\frak S$
such that all natural families of walks on the slit plane with steps
in $\frak S$ have an algebraic \gf .

\medskip
As we did for the cycle lemma result, we can state explicitly the
result on one-dimensional walks that is underlying Propositions~\ref{greg}
and~\ref{propo-cycle-2}. We go back to the framework of 
Lemma~\ref{1D-cycle}.
\begin{Propo}
For $i\ge 1$, let $S_i \equiv S_i(\lambda)$ denote the \gf \ for walks
that start from $0$, end at $i$, and always have a positive level once
they have left their starting point.  This series can be computed
by induction on $i$ via 
the following identity: 
$$\sum_{k =1}^i\frac{(-1)^{k-1}}{k} \sum_{i_1+\cdots+i_k=i\atop i_1
>0, \ldots , i_k >0}
S_{i_1}\cdots S_{i_k} =  [x^{i}] \log \frac 1 {1-\sum_{i \in \frak P} \lambda_i
x^i}.$$
\end{Propo}

\bigskip We still have one thing to do in this section: explain its
title. 
The  way I discovered the first identity of
Proposition~\ref{greg} is not the one that is presented above. The
short and nice proof presented above was ``observed'' by
Gilles Schaeffer, probably less than two minutes after I showed him
the identity...
 The original proof was  longer, and was inspired 
by a probabilistic argument due to Greg Lawler~\cite{lawler},
which we  describe here in combinatorial terms. 

\begin{figure}[ht]
\begin{center}
\epsfig{file=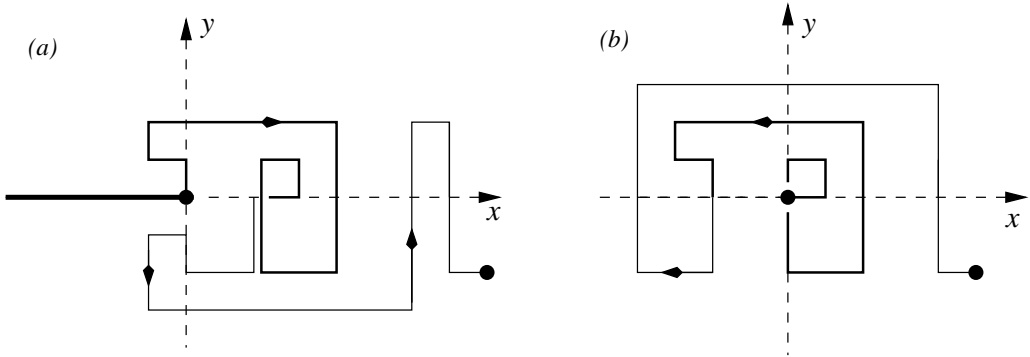}
\end{center}
\caption{Factorization of words of ${\cal S}$ and $\cal W$.}
\label{factor}
\end{figure}

Let $w$ be a non-empty walk on the slit plane, on the ordinary square
lattice. Let  
$j = \min \{i\ge 0 : w \hbox{ visits } (i,0)\}.$
If $j=0$, then $w$ completely avoids the $x$-axis. Otherwise, we factor $w$ at
the last time it visits $(j,0)$ (see Figure~\ref{factor}.$a$). This gives
\beq \GS  = \GA + \GF \ \GS , \label{factorS}\eeq
where 
$\cal A$ is the set of walks that never return to
the $x$-axis and 
$\cal F$ is the set of words $w$ satisfying:
$$|w|_n=|w|_s, \ \delta(w) >0, \  \hbox{ and for all factorizations } w=uv \hbox{ such
that } u \not = \epsilon, \hbox{ either } |u|_n \not = |u|_s \hbox{ or
} \delta(v) \le 0.$$
Let us now consider the set $\cal W$  of all
square lattice walks. 
A walk $w \in {\cal W}$  is of one of the following types:
\begin{itemize}
\item either $w$ avoids the half-line ${\cal H}'=\{(k,0): k<0 \}$; by looking
at the last time when $w$ visits $(0,0)$, we observe that the \gf \
of such walks is ${\GL} \ {\GS}$, where $\GL$ counts loops,
\item or   $w$ meets the half-line ${\cal H}'$; in this case, we
factor $w$ at the first time it meets  ${\cal H}'$ (Figure~\ref{factor}.$b$).
\end{itemize}
This gives
\beq \GW = \GL \ \GS + \bar \GF \ \GW\label{factorW}\eeq
 where $\overline {\cal F}$ is  the
language of words $w$ satisfying
$$|w|_n=|w|_s, \ \delta(w)<0,  \ \hbox{ and for all factorizations } w=uv
\hbox{ such 
that } v \not = \epsilon, \hbox{ either } |u|_n \not = |u|_s \hbox{ or
}\delta(u) \ge 0.$$

It is obvious on  Figure~\ref{factor} that the walks of
$\overline {\cal F} $ are obtained by reversing the direction in walks of $\cal
F$. In terms of words, the words of $\overline {\cal F} $ are obtained
by reading the words of $\cal F$ from right to left, replacing each
occurrence of $o$ (resp. $n,e,s$) by a letter $e$ (resp. $s,o,n$).
Naturally, this can also be checked on the definitions of the sets
$\cal F$ and $\overline{\cal F}$.
Let us convert~\Ref{factorS}
and~\Ref{factorW} into identities on length \gfs , then eliminate
$F$. 
We obtain:
 $$ L(t)S(1,1;t)^2 = A(1,1;t) W(1,1;t).
$$
But the set $\GA$ of walks avoiding the $x$-axis is simply related to
the set $\GB$ of bilateral walks: indeed, by looking at the last time
an ordinary walk visits the $x$-axis, we obtain $\GW = \GB \ \GA$. 
Hence, by Lemma~\ref{lemmeMB}, the above identity reads $
L(t)S(1,1;t)^2 =(1-4t)^{-3/2}$.  This identity, formulated
in probabilistic terms, is Proposition~2.4.6 in~\cite{lawler}. It is
obviously not sufficient to characterize $S(1,1;t)$, but it {\em is\/}
sufficient to prove that the number of $n$-step walks on the slit
plane grows like $4^n n^{-1/4}$, and that's what Lawler does.

 Our first variation on Lawler's approach is to take into account
the endpoint of the walks;  more variables should
give more information. Using~\Ref{factorS}
and~\Ref{factorW}, we thus obtain:
 $$ S(\bx, \by;t)L(t)S(x,y;t) = A(\bx,\by;t) W(x,y;t).$$
This identity  follows from Proposition~\ref{greg}, given that $W=AB=1/K$. 
 Our second variation comes from the
observation that it is difficult to work with series that, like
$S(x,y;t)$, contain both positive and negative powers of $x$ (and
$y$). That is why we focus on the set ${\cal S}_0$ 
  of bilateral walks on the slit plane: their
final abscissa is always positive, and their final ordinate is of
course $0$. The restriction of Lawler's argument to bilateral walks
is straightforward. It boils down to imposing this condition to all
walks under consideration, thus replacing $\cal S$ by  ${\cal S}_0$,
$\cal A$ by $\{\epsilon\}$, and $\cal W$ by $\cal B$, the set of
bilateral walks. This is how we obtained the first identity of
Proposition~\ref{greg}. 
Our third and last  observation is that this can be
done with any set of steps.

\section{The complete solution}
\label{section-complete}

The factorization of $B(x;t)$ given in Proposition~\ref{greg},
obtained by a combinatorial argument, justifies our interest in the following
factorization lemma. This lemma also played a key role in~\cite{prep-GS},
for reasons that were far from being as clear as here.
\begin{Lemma}[The factorization lemma~\cite{prep-GS}]
\label{factorisation-lemma}
Let $B(x;t)$ be a series in $t$  with coefficients in $\rs
[x, \bx]$, and assume  $B(x;0)=1$. There exists a unique triple
$(B_0(t),B_+(x;t), B_-(\bx;t))\equiv(B_0,B_+(x), B_-(\bx))$ of 
formal power series in $t$ satisfying the following conditions: 
\begin{itemize}
\item $B(x)=B_0B_+(x) B_-(\bx),$
\item the coefficients of $B_0$ belong to $\rs$,
\item the coefficients of $B_+(x)$  belong to $\rs[x]$,
\item the coefficients of $B_-(\bx)$  belong to $\rs[\bx]$,
\item  $B_0(0)= B_+(x;0)=B_-(\bx;0)=B_+(0;t)=B_-(0;t)=1$.
\end{itemize}
\end{Lemma}
{\bf Proof.}
 Let us take logarithms in the factorization of $B$:
$$ \log B_0(t) + \log B_+(x; t) + \log B_-(\bx;t) =\log
B(x;t). $$
The conditions  $B_+(0;t)=B_+(x;0)=1$ force the series $\log
 B_+(x;t)$ to be a multiple of $x$. Similarly, $\log B_-(\bx;t)$ must
 be a multiple of $\bx$.
But $\log B(x;t)$ can be written  {\em in a unique  way\/}  as
\beq \log B(x;t)= L_0(t) + xL_1(x;t) + \bx L_2(\bx;t)\label{log+-} \eeq
where $L_0$ does not depend on $x$,  $L_1$ is a series in $t$ with
{polynomial}  coefficients in $x$, and $L_2$ is a series in $t$ with
{polynomial}  coefficients in $\bx$. This forces 
\beq \log B_0(t) =
L_0(t), \quad  \log B_+(x;t)=xL_1(x;t) \quad  \hbox{and } \quad 
\log B_-(\bx;t)=\bx L_2(\bx;t).\label{solution-logs}\eeq
  This shows the existence and
uniqueness of the factorization.
\cqfd
This lemma implies that the first equation of Proposition~\ref{greg}
characterizes completely the series $S_0(x;t)$ and $L(t)$, and reduces
the solution of any slit plane model to 
the explicit  factorization of an algebraic series.
\begin{Theorem}[The complete solution]
\label{theoreme-general}
Let $\frak S$ be  an arbitrary set of steps. Let $B(x;t)$ be the \gf \
for bilateral walks with steps in $\frak S$. The \gf \ $S(x,y;t)$,  $S_0(x;t)$
and $L(t)$, which count various families of walks on the slit plane
(general walks, bilateral walks and loops,
respectively) can be expressed in terms of the polynomial
$K(x,y;t)$ (given by~{\em \Ref{kernel}}) and the canonical factorization
of $B(x;t)$: 
\begin{eqnarray*}
L(t)&=&B_0(t),\\
S_0(x;t)&=& B_+(x;t),\\
S(x,y;t)&=&\displaystyle \frac{1}{K(x,y;t) B_0(t) B_-(\bx;t)}.
\end{eqnarray*}
The logarithm of each of these series is D-finite.
\end{Theorem}
The holonomy results follows from the first point of
Proposition~\ref{nature-fact} below. More specifically,
 we shall be concerned in the sequel of the paper with the {\em algebraic\/}
nature  of the three series given by this theorem. We shall exhibit
some sets of steps $\frak S$ such that all of them are algebraic
(Section~\ref{section-algebraicity}), and some other sets of steps
such that none of them are algebraic, nor even D-finite
(Section~\ref{section-transcendental}). Hence we need to clarify
the nature of the canonical factors of an algebraic series $B(x;t)$.
The following proposition gives a partial answer.
\begin{Propo}[The nature of the canonical factors]\label{nature-fact}
Let $B(x;t)$ be an algebraic series   in $t$  with coefficients in $\rs
[x, \bx]$, such that  $B(x;0)=1$. Let $(B_0(t),B_+(x;t), B_-(\bx;t))$
be its canonical factorization. 
\begin{description}
\item $1.$ The series $\log B_0(t)$, $\log B_+(x;t)$ and $\log B_-(\bx;t)$
are D-finite in all their variables.
\item  $2$. The series $B_0(t)$ is D-finite if and only if the
series 
$$[x^0]\frac{B'(x;t)}{B(x;t)}$$ 
is algebraic (the derivative is taken with respect to $t$).
\item   $3$. The series $B_+(x;t)$ is D-finite in $t$ if and only if the
series $$[x^>]\frac{B'(x;t)}{B(x;t)}:= \sum_{i \ge 1} x^i
[x^i]\frac{B'(x;t)}{B(x;t)}$$ 
is algebraic. In this
case, both  $ B_-(\bx;t)$ and   $B_0(t)$ are also D-finite in $t$.
\item   $4$. If $B(x;t)$ is a polynomial of $\rs[x,\bx,t]$, then $B_0(t)$,
$B_+(x;t)$ and $ B_-(\bx;t)$ are algebraic.
Moreover,
$B_+(x;t)$ is a polynomial in $x$, and $B_-(\bx;t)$ is a polynomial in
$\bx$ (with coefficients in $\rs[[t]]$).
%
\end{description}
\end{Propo}
The proof of this proposition requires the following
result.
\begin{Propo}[Singer~\cite{singer}]\label{singer}
Let $Y(t)$ be power series in $t$ with coefficients in a field
$\GK$ of
characteristic zero. Assume $Y'(t)/Y(t)$ is holonomic. Then $Y(t)$ is
holonomic  if and only if  $Y'(t)/Y(t)$ is actually algebraic.
\end{Propo}
It is easy to prove that the algebraicity of $Y'(t)/Y(t)$ implies the
holonomy of $Y(t)$, but the converse is significantly more difficult.

\medskip
\noindent
{\bf Proof of Proposition~\ref{nature-fact}}\\
\noindent
1. As in the proofs of Lemma~\ref{B-general} and
Proposition~\ref{propo-cycle-1}, we take a new 
indeterminate $\bx$, independent of $x$, and form the \gf \ $\overline
B(x,\bx;t)$ of bilateral walks obtained by giving a weight $x^i t$ to
any step $(i,j)$ such that $i\ge 0$ and  a weight $\bx^i t$ to
any step $(-i,j)$ such that $i< 0$. This series is the diagonal of a
rational series and is algebraic over $\qs(x,\bx,t)$. As the logarithm
of an algebraic series is D-finite,
$\log \overline 
B(x,\bx;t)$ is D-finite in its three variables. By
Eqs.~(\ref{log+-}--\ref{solution-logs}), the 
series $\log B_0(t)$ is the diagonal of $\log \overline
B(x,\bx;t)$  in $x$ and $\bx$, evaluated at $x=1$. Similarly, $\log
B_+(x;t)$ is the diagonal  of $u/(1-u) \times \log \overline
B(x,u\bx;t)$ in $x$ and $u$, evaluated at $\bx=1/x$. Finally,  $\log B_-(\bx;t)$ is the diagonal  of $u/(1-u) \times \log \overline
B(ux,\bx;t)$ in $\bx$ and $u$, evaluated at $\bx=1/x$. The properties
recalled in 
Section~\ref{rappels-DF} imply that these three series are D-finite in
all their variables. 

\medskip
\noindent 2. By differentiating the first identity
of~\Ref{solution-logs} with respect to $t$, we obtain: 
$$\frac{d}{dt} \log B_0(t)= \frac{B'_0(t)}{B_0(t)}=\frac{d}{dt} [x^0]
\log B(x;t)= 
 [x^0]\frac{\partial }{\partial t}\log B(x;t) =
[x^0]\frac{B'(x;t)}{B(x;t)}.$$ 
This series is clearly D-finite (it is the derivative of a D-finite
series). The announced result follows from Proposition~\ref{singer},
applied to $Y=B_0$ and $\GK= \cs$.

\medskip
\noindent 3. Similarly, 
$$\frac{B'_+(x;t)}{B_+(x;t)}=[x^>]\frac{B'(x;t)}{B(x;t)},$$
and  Proposition~\ref{singer}, applied to
$Y=B_+(x;t)$ and $\GK= \cs(x)$, proves the first statement of point 3.
Now with obvious notations,
$$\frac{B'(x;t)}{B(x;t)}= [x^\le]\frac{B'(x;t)}{B(x;t)}+
[x^>]\frac{B'(x;t)}{B(x;t)},$$
and $B'/B$ is algebraic. Hence, if $[x^>]{B'(x;t)}/{B(x;t)}$ is
algebraic, then so is $[x^\le]{B'(x;t)}/{B(x;t)}$. This series is
a power series in $t$ with coefficients in $\rs[\bx]$.  Setting $\bx
=0$ in this series yields $[x^0]{B'(x;t)}/{B(x;t)}$, which is,
consequently, also algebraic. By difference,
$[x^<]{B'(x;t)}/{B(x;t)}$ is also algebraic; finally
Proposition~\ref{singer}, used in the easy direction, implies that
both $B_0(t)$ and $B_-(\bx;t)$ are D-finite in $t$.

\medskip
\noindent 4. This point was already proved in~\cite{prep-GS}. We
repeat the proof, since we shall use it to solve explicit examples.
 Let the smallest
 exponent of $x$ occurring in $B(x;t)$ be  $-m$. Then $P(x;t)=x^mB(x;t)
$  is a
 polynomial in $x$ and $t$ such that $P(x;0)=x^m$. As a polynomial in
$x$, $P$ has degree, say, $d$, and hence admits $d$ roots,
denoted  $X_1, \ldots , X_d$, which
 are Puiseux series in $t$ with complex coefficients; in particular,
 there exists an integer 
$n\ge 1$ such that all these roots can be written as Laurent series
in the variable    $z=t^{1/n}$ (see~\cite[Theorem~6.1.5]{stanley}). 
Assume that  exactly $k$ of these roots, say
 $X_1, \ldots , X_k$, are finite at $z=0$.
The other $d-k$ roots contain terms of the form $z^{-i}$, with $i>0$.
The polynomial  $P(x;t)$  can be factored as
$$P(x;t)= x^mB(x;t) = \tilde B_0(z) \prod_{i=1}^k (x- X_i)  \prod_{i=k+1}^d
\left(1-\frac{x}{  X_i}\right),$$
where $\tilde B_0(z)$ is, as  the $X_i$, an algebraic function of $z$. For $i>k$,
the series $1/X_i$ equals $0$ at $z=0$. Hence the condition
$P(x;0)=x^m$ implies that $k=m$, that $\tilde B_0(0)=1$, and that the finite
roots  $X_1, \ldots , X_m$ equal $0$ when $z=0$. Let
$$\tilde B_+(x;z)=\prod_{i=m+1}^d \left(1-\frac{x}{
X_i}\right) \ \ \ \hbox{and}\ \ \ 
 \tilde B_-(\bx;z)=\prod_{i=1}^m (1-\bx X_i).$$ Then   the  series
 $\tilde B_0(z), \tilde B_+(x;z)$ and $  \tilde B_-(\bx;z)$ clearly satisfy
all the conditions defining the canonical factorization of $B(x;t)$,
but two: we still 
need to prove that they are actually 
series in $t$ (and not only in  $z=t^{1/n}$) with real (rather than
complex) coefficients. As in  the proof of
Lemma~\ref{factorisation-lemma}, we argue that 
the series  $\log B(x;z^n)= \log B(x;t)$ has a unique decomposition as a sum of
three series (one with positive powers of $z$, one with negative powers
of $z$, one with terms independent of $z$ --- see Eq.~\Ref{log+-}) to prove
that the series  $\tilde B_0(z), \tilde B_+(x;z)$ and $  \tilde
B_-(\bx;z)$ are indeed 
the canonical factors of $B(x;t)$.
\cqfd

\medskip

\noindent{\bf Remark: effectiveness of the solution}

\noindent
 One might wonder how effective the solution of the slit plane models given by
   Theorem~\ref{theoreme-general}  is. We claim that for any given set
   $\frak S$, one can compute algorithmically a linear differential equation
   (in $t$) satisfied by each of the series $\log L(t)$, $\log
   S_0(x;t)$ and $\log    S(x,y;t)$. The reason for this is that all
   the results about rational, algebraic and D-finite series that we
   have used to conclude 
   that $\log S(x,y;t)$ and its siblings are D-finite are  effective,
   and actually implemented in the {\sc Maple} packages {\sc Gfun} and
   {\sc Mgfun}. 

Let us give a few details. One starts from the 
   polynomial $K(x,y;t)$ given by~\Ref{kernel}. By
   converting 
   $1/K$ into partial fractions of $y$, we obtain an expression of
   $B(x;t)$ in terms of the solutions $Y_i$ of $K(x,Y_i;t)=0$. By
   elimination of the $Y_i$, we then  obtain a polynomial equation
   satisfied by the \gf \ 
   $B(x;t)$ for bilateral walks.  We can just as well compute a
   polynomial equation for the three-variable series $\overline
   B(x,\bx;t)$ used in the proof of
   Proposition~\ref{nature-fact}. That is the first step. An algorithm
   that derives from the polynomial equation satisfied by $\overline B$ a
   differential equation (in any of the variables $x, \bx$ or $t$)
   satisfied by $\log \overline B$ is implemented in  {\sc Gfun}. The
   extraction  
   of diagonals is  implemented in  {\sc Mgfun}. One thus obtains
   differential equations for the series $\log B_0$, $\log B_+(x;t)$
   and  $\log B_-(\bx;t)$, and hence for $\log L(t)$ and $\log
   S_0(x;t)$.  Finally, we need an effective version of the closure of
   D-finite series under the sum to obtain a differential equation (in
   $t$, say)    satisfied by $\log S(x,y;t)$: this is provided by the package  {\sc Gfun}. 

This holds for any set $\frak S$.  Moreover, we shall study in
   Section~\ref{section-algebraicity} a family of sets $\frak S$
   such that the \gf \ of bilateral walks is of the form
   $B(x;t)=P(x;t)^{-1/2}$. The last point of
   Proposition~\ref{nature-fact} tells us how to obtain {\em
   explicit expressions\/} of the canonical factors of $B(x;t)$ in
   terms of the roots $X_i$ of $P$.


\section{A step-by-step approach}
\label{section-functional}
Theorem~\ref{theoreme-general} relates three \gfs \ for walks on the
slit plane to the canonical factorization of the series $B(x;t)$ that
counts bilateral walks. In this section, we present another derivation
of this theorem. 
The starting point of this alternative proof is a functional equation
satisfied by the \gf \ $S(x,y;t)$ that counts all walks on the slit
plane. This equation simply translates the fact that a walk  on the slit
plane is obtained by... adding a step to another walk  on the slit
plane. This extremely elementary argument, and the corresponding
equation, were also the starting point of~\cite{prep-GS}. But in that
paper, we were only able to solve the functional equation when the set
$\frak S$ of steps satisfies two simple conditions.
Here, we shall solve the equation by a new and simpler method, which
works for any set of steps $\frak S$. This approach is actually so
simple that I find it somewhat mystifying... It seems that we solve an
a priori non-trivial problem without any clever combinatorial or
algebraic argument. In more positive words, it shows that the crucial
ingredient of the solution is the factorization lemma. 
\begin{Lemma} \label{functional}
Let $S(x,y;t)$ and $L(t)$ denote the \gfs \
for walks  on the slit plane and for loops,  respectively. There exists a
series $\Omega(\bx;t)$ in $\ns[\bx][[t]]$ such that $\Omega(\bx;0)=0$
and
$$K(x,y;t)S(x,y;t) =1-\Omega(\bx;t), \quad 
L(t)=\frac{1}{1-\Omega(0;t)}$$
with, as above,
\beq K(x,y;t)=1-t\sum _{(i,j) \in \frak S} x^i y^j . \label{kernel-general}\eeq
\end{Lemma}
{\bf Proof.}
We obtain an equation for the series $S(x,y;t)$ by saying
that a walk of length $n$  is
obtained by adding a step to another walk of length
$n-1$. However,  this procedure sometimes
produces a {\em bridge\/},  that is, a nonempty walk that 
starts  at
$(0,0)$, ends on the half-line $\cal 
H$, but otherwise avoids $\cal H$. Hence, denoting by $\Omega(\bx;t)$
the \gf \ for bridges, we have:
$$S(x,y;t)=1+tS(x,y;t) \left(\sum _{(i,j) \in \frak S} x^i y^j \right)
-\Omega(\bx;t),$$
which is exactly the announced equation for $S(x,y;t)$.
Now a loop  is a sequence of bridges, all ending at $(0,0)$. The
equation for $L(t)$ follows.
\cqfd
\begin{Coro}
Theorem~{\em \ref{theoreme-general}} holds.
\end{Coro}
{\bf Proof.} We have established
\beq S(x,y;t)= \frac{1-\Omega(\bx;t)}{K(x,y;t)}. \label{eq1}\eeq
Note that $1/K(x,y;t)$ counts all walks with steps in $\frak S$. 
Extracting the coefficient of $y^0$ in this equation gives:
$$S_0(x;t)= B(x;t)({1-\Omega(\bx;t)}),$$
where $B(x;t)$ is the \gf \ for bilateral walks. We claim that this
immediately implies Theorem~\ref{theoreme-general}! Indeed, the
uniqueness  of the 
canonical factorization of $B(x;t)$ (Lemma~\ref{factorisation-lemma}),
combined with  the properties of
$S_0(x;t)$ and $\Omega(\bx;t)$, 
forces
$$ S_0(x;t)= B_+(x;t) $$
(one of the three results of Theorem~\ref{theoreme-general}), {and}
\beq 1-\Omega(\bx;t)=\frac{1}{B_0(t)B_-(\bx;t)}. \label{eq2}\eeq
The latter identity, plugged in~\Ref{eq1}, yields
$$ S(x,y;t)= \frac{1}{K(x,y;t)B_0(t)B_-(\bx;t)}$$
which is the third statement of
Theorem~\ref{theoreme-general}. Finally, \Ref{eq2} and
Lemma~\ref{functional} imply that $ L(t)=B_0(t),$ which completes this
new proof of Theorem~\ref{theoreme-general}.
\cqfd

\noindent{\bf  Remark: a connection with Spitzer's book}\\
 It is well-known that {\em everything\/} is in Spitzer's book, {\em
   Principles of random walks\/}~\cite{spitzer}. The results of
   Sections~\ref{section-complete} and~\ref{section-functional} are, to a certain extent,   not an exception
   to this rule. In Chapter $4$ of his book, Spitzer studies
   generic one-dimensional random walks $0=S_0, S_1, S_2\ldots$ where
   $S_i \in\zs$ for all $i$. He is interested in the following
   stopping times:
$$T=\min\{ 1\le n \le \infty : S_n > 0\} \quad \hbox{and} \quad
T'=\min\{ 1\le n \le \infty : S_n \ge 0\}.$$
Our two-dimensional
combinatorial problem can be turned into a probabilistic one as follows:
assume each step of $\frak S$ is taken
with probability $1/k$, where $k$ is the cardinality of $\frak
S$. Then, if we
   only look at the positions where this two-dimensional random walk
hits the $x$-axis, we effectively obtain a 
   one-dimensional random walk with steps in $\zs$. 
The series $f_e(t,x)$, $f_i(t,x)$ and
   $c(t)$ defined in Definition~17.D3 of~\cite{spitzer} (and also
used in~\cite{kesten}) are
    then the counterparts of $1/B_+(x;t), 1/B_-(\bx;t)$ and
   $1/B_0(t)$, respectively. Proposition~17.P4 in~\cite{spitzer} is the
   counterpart of $B=B_- B_+ B_0$. The first two results of
   Proposition~17.P5 are related to the enumeration of bridges; more
   precisely, Eq.~17.P5(b) is the counterpart of Eq.~\Ref{eq2} above.
The next two results of
   Proposition~17.8 are related to the enumeration of walks avoiding
   $\cal H$; more
   precisely, Eq.~17.P5(d) is the counterpart of the second identity of
   Theorem~\ref{theoreme-general}. And so on! Of course, there are some
   differences between our treatment and Spitzer's: we deal with exact
   enumeration rather than probabilities, and a unit time in our model
   does not coincide with a unit time in Spitzer's problem. But these
   differences are minor, and one has to realize that we have (as
   always?) solved our two-dimensional model by reducing it to a
   one-dimensional question. What happens between two visits to the
   $x$-axis has no importance. However, what we definitely add to
   Spitzer's treatment is that our results are effective; the next two
   sections will provide many explicit examples.

\section{Walks with small height variation: algebraicity}
\label{section-algebraicity}
In this section, we prove that if the set of steps  $\frak S$ satisfies
the following {\em small height variation\/} condition:
 $$ \hbox{for all } (i,j) \in \frak S, \quad |j|\le 1, $$
then all the \gfs \ for walks on the slit plane defined in
Section~\ref{definitions} are algebraic. This condition prevents a
walk from crossing the half-line $\cal H$ without hitting it. 
Three examples will be solved explicitly: the 
ordinary square lattice and the diagonal square lattice studied in
Section~\ref{section-cycle}, and the triangular lattice of
Figure~\ref{chemin-triangle}.   

Let $\frak S$ be a set of steps satisfying the small height variation
condition. We define three Laurent polynomials in $x$, denoted $
A_{-1}(x)$, $A_0(x)$ and $A_1(x)$, by
$$ A_j(x)= \sum_{(i,j) \in \frak S} x^i.$$
The polynomial $K(x,y;t)$, given
by~\Ref{kernel-general}, can be rewritten 
$$ K(x,y)=1-t\by A_{-1}(x) - tA_0(x)-ty A_1(x).$$

\subsection{General results}
The following result was proved in~\cite{prep-GS} under the symmetry
hypothesis $A_1(x)=A_{-1}(x)$. The new approaches developped in this
paper  show that this assumption is unnecessary. 
\begin{Theorem}\label{complete-small}
Let $\frak S$ be a set of steps with small height variations. Let
$\Delta(x;t)$ be the following polynomial in $x, \bx$ and $t$:
\beq  \Delta(x;t)= (1-tA_0(x))^2 -4t^2A_1(x)A_{-1}(x).\label{Delta-def} \eeq
Then the \gf \ for bilateral walks is
$$B(x;t)= \frac{1}{\sqrt {\Delta(x;t)}}.$$ 
Let $(\Delta_0(t),\Delta_+(x;t), \bar \Delta_-(\bx;t))$ be the
canonical factorisation of
$\Delta(x;t)$, defined  in  Lemma~{\em \ref{factorisation-lemma}}. Then
the \gf \ for 
walks on the slit plane with steps in $\frak S$ is
$$S(x,y;t)= \frac{\displaystyle \sqrt{\Delta_0(t)
\Delta_-(\bx;t)}}{K(x,y;t)}.$$ 
The \gf \ for bilateral walks on the slit plane is
\beq S_0(x;t) =
\frac{1}{\sqrt {\Delta_+(x;t)}}.\label{S0-small}\eeq
The \gf \ for loops is
$$L(t)=\frac{1}{\sqrt {\Delta_0(t)}}.$$ 
All  these series  are algebraic. 
\end{Theorem}
{\bf Proof.}
We proceed as  in
Lemma~\ref{lemmaMBdiag} to enumerate bilateral walks. 
The \gf \ $M(x;t)$ for bilateral walks that
stay on or above the $x$-axis satisfies
\begin{eqnarray*}
M(x;t)&=&1+tA_0(x)M(x;t)+
tA_1(x)M(x;t)tA_{-1}(x)M(x;t)\\
&=&1+tA_0(x)M(x;t)+t^2A_1(x)A_{-1}(x)M(x;t)^2.
\end{eqnarray*}
 The  \gf \ $B$ for bilateral walks is related to $M$ by
$$B(x;t)=1+tA_0(x)B(x;t)+2t^2A_1(x)A_{-1}(x)M(x;t)B(x;t).$$
(The factor $2$ comes from the fact that primitive bilateral words
starting with a down step are obtained by reading from right to left
primitive bilateral words starting with an up step.) The expression of
$B(x;t)$ follows.
 Theorem~\ref{theoreme-general} provides the values of
$S(x,y;t)$, $S_0(x;t)$ and $L(t)$. By Proposition~\ref{nature-fact},
the series $\Delta_0(t)$,  
$\Delta_+(x;t)$ and $\Delta_-(\bx;t)$ are algebraic, and so are
$S(x,y;t)$, $S_0(x;t)$ and $L(t)$. 
\cqfd

As in~\cite{prep-GS}, we derive from this central result that the \gf
\ for walks on the slit plane ending at a prescribed position $(i,j)$
is also algebraic. This, we believe, is a surprising result, as it is
not true for ordinary walks on the square lattice, nor for walks
avoiding the whole $x$-axis.
\begin{Theorem}\label{Sij-small}
  For all $j \in \zs$, the \gf \
$S_j(x;t)$ for walks on the slit plane ending at ordinate $j$ is
algebraic. More precisely,  for $j\ge 0$,
$$S_j(x;t)= \frac{1}{\sqrt{\Delta_+(x;t)}}
\left(\frac{Z(x;t)}{A_{-1}(x)}\right)^j,$$
$$S_{-j}(x;t)= \frac{1}{\sqrt{\Delta_+(x;t)}}
\left(\frac{Z(x;t)}{A_{1}(x)}\right)^j,$$
where
$$Z(x;t)=\frac{1-tA_0(x)-\sqrt{\Delta(x;t)}}{2t}.$$
and $\Delta(x;t)$ is given by~{\em \Ref{Delta-def}}.
 The \gf \ $S_{i,j}(t)$ for walks on the slit plane ending at $(i,j)$
is also algebraic for any $(i,j)$:
\begin{description}
\item -- If $j= 0$, Eq.~{\em\Ref{S0-small}}
shows that  $S_{i,j}(t)$  belongs to
the extension of $\qs(t)$ generated by  the
coefficients of $\Delta_+(x;t)$  (seen as
a polynomial in $x$).
\item
-- For $j > 0$, it belongs to
the extension of $\qs(t)$ generated by $\sqrt{\Delta_0(t)}$, the
coefficients of $\Delta_+(x;t)$ and  $\Delta_-(x;t)$ (seen as
polynomials in $x$),  and
also the series $\sqrt{\Delta_+(\alpha _i;t)}$, where the numbers $\alpha _i$
are the roots of the polynomial $A_{-1}(x)$, and the algebraic numbers
$\alpha _i$ themselves. 
\item
-- A similar statement holds if $j<0$, upon
replacing  $A_{-1}(x)$ by  $A_{1}(x)$.
\end{description}
\end{Theorem}
{\bf Proof.} The proof is very similar to what we did
in~\cite{prep-GS}. The series $S_j(x;t)$ is the coefficient of $y^j$
in $S(x,y;t)$. To extract this coefficient from the expression of
$S(x,y;t)$ given in Theorem~\ref{complete-small}, we need to convert
the rational 
function $1/K(x,y;t)$ into partial fractions of $y$. As a polynomial
in $y$, $K(x,y;t)$ has two roots. One of them is $Z(x;t)/A_{1}(x)$,
and the other is $A_{-1}(x)/Z(x;t)$. Hence
\begin{eqnarray*}
\frac{1}{K(x,y;t)}&=& \frac{y
Z(x;t)}{t(yA_1(x)-Z(x;t))(A_{-1}(x)-yZ(x;t))}\\
&=& \frac{1}{\sqrt{\Delta(x;t)}}\left(\frac{1}{1-\by \frac{Z(x;t)}{A_1(x)}}+ \frac{1}{1-y \frac{Z(x;t)}{A_{-1}(x)}}-1 \right).
\end{eqnarray*}
As $Z(x;t)$ is a series in $t$ (with coefficients
in $\rs[x,\bx]$) such that $Z(x;0)=0$, 
each term in the above expression is a well-defined series in $t$, and
it is now easy to extract the coefficient of $y^j$ from $1/K(x,y;t)$,
and hence from $S(x,y;t)$. The announced
expressions of $S_j$ and $S_{-j}$ follow. From this point, the
derivation of $S_{i,j}(t)$ exactly copies what was done
in~\cite[Section~3.2]{prep-GS}, and we do not  repeat
the argument.
\cqfd

\subsection{Three examples}\label{alg-examples}
The first two examples studied below (the ordinary square lattice and
the diagonal square lattice) were already solved in~\cite{prep-GS}.
\begin{Theorem}[The ordinary square lattice]\label{dim2}
The \gf \ $S(x,y;t)$ for walks on the slit plane,  the
\gf \ $S_0(x;t)$ for bilateral walks on the slit plane
and the \gf \ $L(t)$ for loops are given by:
$$
S(x,y;t) =\frac{\left(1-2t(1+\bar x)+\sqrt{1-4t}\right)^{1/2}
                   \left(1+2t(1-\bar x)+\sqrt{1+4t}\right)^{1/2}}
             {2(1-t(x+\bar x+y+\bar y))},
$$
$$S_0(x;t)=  \frac{2t}{ \left(2t-x(1-2t-\sqrt{1-4t})\right)^{1/2}
\left(2t-x(1+2t-\sqrt{1+4t})\right)^{1/2}},$$
$$L(t)= \frac{(\sqrt{1+4t}-1)(1-\sqrt{1-4t})}{4t^2}= \sum_{n \ge 0}
(2.4^nC_n -C_{2n+1})t^{2n},$$
where $C_n = {{2n} \choose n}/(n+1)$ is the $n$th Catalan number.
These series are algebraic of degree $8$, $8$ and $4$ respectively. 

 For any $(i,j)$, the \gf \ $S_{i,j}(t)$ for walks on the slit
plane ending at $(i,j)$ belongs to $\qs(t,\sqrt{1-4t}, \sqrt{1+4t})$.
\end{Theorem}
{\bf Proof.}
We apply Theorem~\ref{complete-small} with $A_0(x)=x+\bx$ and
$A_1(x)=A_{-1}(x)=1$. Hence $ {\Delta(x;t)}=\left( 1-t(x+\bx
+2)\right )\left( 1-t(x+\bx 
-2)\right )$. In order to compute the canonical factorization of  $
{\Delta}$, we apply the procedure described in the proof of 
Proposition~\ref{nature-fact}.  The
polynomial $\Delta$ has four roots $X_i$, $1 \le i \le 4$, which are
quadratic functions of $t$ and can be computed explicitly. Let $C(t)$
denote the \gf \ for Catalan numbers:
\beq C(t)=\frac{1-\sqrt{1-4t}}{2t}=\sum_{n \ge 0}\frac{1}{n+1}{2n \choose
n} t^n.\label{catalan} \eeq
Then $X_1= C(t)-1$ and $X_2=1-C(-t)$ are the two roots that are finite
at $t=0$, and by symmetry of $\Delta$ in $x$ and $\bx$, the two other
roots are 
\beq X_3= (C(t)-1)^{-1} = \frac{1-2t+\sqrt{1-4t}}{2t} 
\quad \hbox{and  } \quad  X_4=(1-C(-t))^{-1}=
\frac{1+2t+\sqrt{1+4t}}{2t}.\label{X34}\eeq
Hence the canonical factorization of $\Delta$ is such that
\begin{eqnarray}
  \Delta_+(x;t)=\Delta_-(x;t)&=&(1-x X_3^{-1})(1-x X_4^{-1})
\label{Delta-formel} \\ 
&=&\left( 1-x (C(t)-1)\right) \left(
1+x(C(-t)-1) \right).
\label{Delta-def1}
\end{eqnarray}
Taking the coefficient of $x^2$ in the relation
$\Delta(x)=\Delta_0(t)\Delta_+(x;t) \Delta_-(\bx;t)$ yields
\begin{eqnarray}
 \Delta_0(t)&=& t^2 X_3 X_4 \label{D-formel}\\
&=&(C(t) C(-t))^{-2},\label{delta-Delta} 
\end{eqnarray}
as $(C(t)-1) = tC(t)^2$.
Now by Theorem~\ref{complete-small}, Eqs.~\Ref{Delta-formel}
and~\Ref{D-formel}, 
$$S(x,y;t) = \frac{\sqrt{t^2  X_3 X_4(1-\bx X_3^{-1})(1-\bx
X_4^{-1})}}
{K(x,y;t)} = \frac{\sqrt{(t X_3-t\bx)(t X_4-t\bx)}}
{1-t(x+\bx+y+\by)}.$$ 
The announced expression of $S$ now follows from~\Ref{X34}. Finally,
Theorem~\ref{complete-small}, combined with the explicit
expressions~\Ref{Delta-def1} and~\Ref{delta-Delta}, gives the
expressions of $S_0(x;t)$ and $L(t)$. 

We now apply Theorem~\ref{Sij-small}: the polynomial
$A_1(x)=1$ having no root, the series $S_{i,j}(t)$ belongs, by~\Ref{Delta-def1}
and~\Ref{delta-Delta}, to $\qs(t, C(t), C(-t))$, that is, to
$\qs(t,\sqrt{1-4t},\sqrt{1+4t})$. 
\cqfd

\begin{Theorem}[The diagonal square lattice]\label{dim2-diag}
The \gf \ $S(x,y;t)$ for walks on the slit plane,  the
\gf \ $S_0(x;t)$ for bilateral walks on the slit plane
and the \gf \ $L(t)$ for loops are given by
$$S(x,y;t)=
\frac{\left(1-8t^2(1+\bx ^2)+\sqrt{1-16t^2}\right)^{1/2}}
{\sqrt{2}(1-t(x+\bx)(y+\by))}.$$ 
$$
S_0(x;t)= \frac{2\sqrt{2}t}{
\left(8t^2-x^2
(1-8t^2-\sqrt{1-16t^2}
)\right)^{1/2}}, $$
and
$$L(t)=\frac{1-\sqrt{1-16t^2}}{8t^2} =\sum_{n \ge 0} 4^n C_n t^{2n},$$ 
where $C_n = {{2n} \choose n}/(n+1)$ is the $n$th Catalan number.
These series are algebraic of degree $4$, $4$ and $2$ respectively.

 For any $(i,j)$, the \gf \ $S_{i,j}(t)$ for walks on the slit
plane ending at $(i,j)$ belongs to $\qs(t,\sqrt{1-4t}, \sqrt{1+4t})$.
\end{Theorem}
{\bf Proof.}
For the diagonal square lattice,  $A_0(x)=0$ and
$A_1(x)=A_{-1}(x)=x+\bx$. Again, the 
polynomial $ {\Delta(x)}= 1-4t^2(x+\bx)^2$   has four roots $X_i$, $1 \le i \le 4$, which are quadratic
functions of $t$ and can be expressed explicitly in terms of the
Catalan \gf \ \Ref{catalan}.  More precisely, $X_{1,2} =\pm 2t C(4t^2)$ and
\beq X_{3,4}=\pm X_1^{-1} = \pm\frac{1+\sqrt{1-16t^2}}{4t} 
.\label{X34-bis} \eeq
We then follow the same steps as above. The formal
expression of $\Delta_+(x;t)$  in
terms of the roots $X_i$ remains unchanged,  but its actual
value is of course different:
\begin{eqnarray} 
\Delta_+(x;t) =\Delta_-(x;t) &=& (1-x X_3^{-1})(1-x
X_4^{-1})\label{Delta-bis-formel}\\ 
&=&{1-4t^2x^2C(4t^2)^2}.
\label{Delta-bis}
\end{eqnarray}
 Taking the coefficient of $x^2$ in the relation $\Delta(x;t)=\Delta_0(t)\Delta_+(x;t)
\Delta_-(\bx;t)$ yields now
\begin{eqnarray}
 \Delta_0(t)&=& -4t^2 X_3 X_4 \label{D-formel-2}\\
&=& C(4t^2)^{-2}.\label{delta-Delta-2} 
\end{eqnarray}
The expression of $S(x,y;t)$ now follows from  Theorem~\ref{complete-small},
Eqs.~\Ref{Delta-bis-formel}, \Ref{D-formel-2}
and~\Ref{X34-bis}.  Finally,
Theorem~\ref{complete-small}, combined with the explicit
expressions~\Ref{Delta-bis} and~\Ref{delta-Delta-2}, gives the
expressions of $S_0(x;t)$ and $L(t)$.

\begin{figure}[ht]
\begin{center}
\epsfig{file=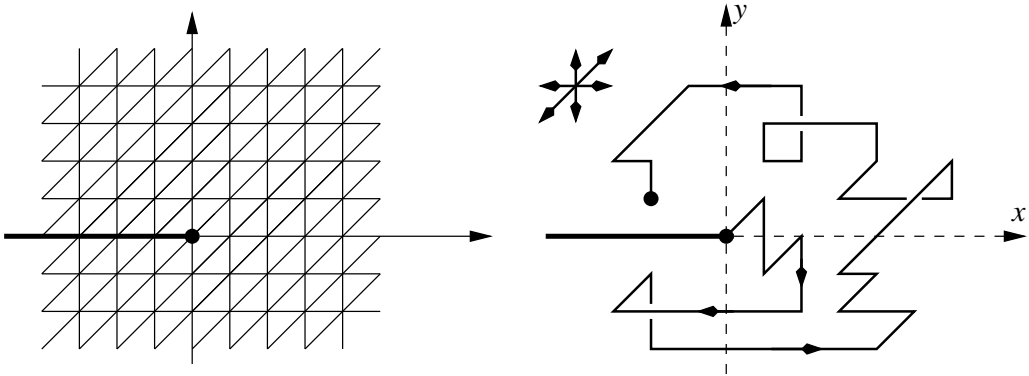}
\end{center}
\caption{A walk on the slit plane with steps
in $\{(0,1),(1,0),(1,1),(0,-1),(-1,0),(-1,-1)\}$.} 
\label{chemin-triangle}
\end{figure}

We now apply Theorem~\ref{Sij-small}. The extension of
$\qs(t)$ generated by $\sqrt{\Delta_0(t)}$ and the coefficients of
$\Delta_+(x;t)$ and $\Delta_-(\bx;t)$ is, by~\Ref{delta-Delta-2}
and~\Ref{Delta-bis}, the field
$\qs(t,C(4t^2))=\qs(t,\sqrt{1-16t^2})$. But the polynomial
$A_1(x)=A_{-1}(x)=x+\bx $ has two  roots at $x=\pm i$, and
$$\sqrt{\Delta_+(\pm i;t)} = \sqrt{1+4t^2 C(4t^2)^2}= \sqrt {C(4t^2)}
=
\frac{\sqrt{1+4t}-\sqrt{1-4t}}{4t} ,$$
so that the series $S_{i,j}(t)$ finally belongs to
$\qs(t,\sqrt{1-4t},\sqrt{1+4t})$. 
\cqfd
We conclude this section by a more complicated example: the triangular
lattice of Figure~\ref{chemin-triangle}. Each 
edge of the lattice can be traversed in both directions, so that the
set of steps is $\frak S=
\{(0,1),(1,0),(1,1),(0,-1),(-1,0),(-1,-1)\}$. Hence $A_0(x)=x+\bx$,
$A_1(x)=1+x$ and $A_{-1}(x)=1+\bx$, so that
$$\Delta(x;t)= (1-t(x+\bx))^2-4t^2(1+x)(1+\bx).$$
This polynomial has  again four roots,
but it is now irreducible, whereas it split into two quadratic terms
in the previous examples. We shall express our results in terms of the
unique series $J(t)=t+ O(t^2)$ satisfying
\beq J=t \ \frac{1-2J+6J^2-2J^3+J^4}{(1-J)^2}=t+5t^3+8t^4+ \cdots
\label{Jeq}\eeq
This series can be expressed with radicals:
\beq J=1+\frac{ 1-2t+\sqrt{(1+2t)(1-6t)}}{4t}\left(1
-\sqrt{\frac{1-\sqrt{(1+2t)(1-6t)}}{2t}}\right).\label{J-expr}\eeq
Observe that 
$$\frac{1-2t-\sqrt{(1+2t)(1-6t)}}{8t^2} = \sum_{n\ge 0} M_n 2^n t^n$$
where $M_n$ is the $n$th {\em Motzkin\/}  number.

\begin{Theorem}[The triangular lattice]
The \gf \ $S(x,y;t)$ for walks on the slit plane,  the
\gf \ $S_0(x;t)$ for bilateral walks on the slit plane
and the \gf \ $L(t)$ for loops are given by
$$
S(x,y;t) = \frac{t\sqrt{(1-J)^2-2J(1+J^2)\bx +J^2(1-J)^2 \bx
^2}}{J(1-J)(1-t(x+\bx+y(1+x)+\by(1+\bx)))} ,
$$
$$
S_0(x;t)= \frac{1-J}{\sqrt{(1-J)^2-2J(1+J^2)\bx +J^2(1-J)^2 \bx ^2}}, $$
and
$$L(t)=\frac{J} t,$$ 
where $J$ is the series in $t$ defined by~{\em \Ref{J-expr}}.
These series are algebraic of degree $8$, $8$ and $4$ respectively.

 For any $(i,j)$, the \gf \ $S_{i,j}(t)$ for walks on the slit
plane ending at $(i,j)$ belongs to $\qs(J)$.
\end{Theorem}
{\bf Proof.} As $\Delta(x;t)$ is quadratic in $x+\bx$, its four roots
can be easily 
expressed in terms of square roots. One finds:
$$X_1=\frac{1+2t+2\sqrt{t(1+3t)}
-\sqrt{(1+2t)(1+6t+4\sqrt{t(1+3t)})}}{2t}= t-2t\sqrt t +2t^2
+\cdots,$$
$$X_2=\frac{1+2t-2\sqrt{t(1+3t)}
-\sqrt{(1+2t)(1+6t-4\sqrt{t(1+3t)})}}{2t}= t+2t\sqrt t +2t^2 +\cdots,$$
$$X_3=\frac{1} {X_1} =\frac{1+2t+2\sqrt{t(1+3t)}
+\sqrt{(1+2t)(1+6t+4\sqrt{t(1+3t)})}}{2t}= \frac 1 t+ \frac 2{\sqrt t
}  +2 +\cdots,$$
$$X_4=\frac{1} {X_2} =\frac{1+2t-2\sqrt{t(1+3t)}
+\sqrt{(1+2t)(1+6t-4\sqrt{t(1+3t)})}}{2t}= \frac 1 t- \frac 2{\sqrt t
}  +2 +\cdots$$
Hence, once again,
\beq
  \Delta_+(x;t)=\Delta_-(x;t)=(1-x X_3^{-1})(1-x X_4^{-1})
=(1-x X_1)(1-x X_2). \label{delta-triangle} \eeq
Taking the coefficient of  $x^2$ in the relation
$\Delta(x;t)=\Delta_0(t)\Delta_+(x;t) 
\Delta_-(\bx;t)$ yields now
 $\Delta_0(t)= t^2 X_3 X_4 $.
By Theorem~\ref{complete-small}, the \gf \ for loops is thus
$L(t)=\sqrt{X_1X_2}/t$. Rather than juggling with the expressions of
the $X_i$'s, we can obtain directly an equation for $L(t)$ by
eliminating $X_1$ and $X_2$ in the equations
$$\Delta(X_1;t)=0, \quad \Delta(X_2;t)=0, \quad t^2L(t)^2=X_1X_2.$$
The resulting equation factors into several polynomials. The fact that
$L(t)=1+5t^2+O(t^3)$ tells us which of these factors cancels $L(t)$, and we
actually end up with $L(t)=J/t$. That is,
\beq X_1X_2=J^2 \quad \hbox{and} \quad
\Delta_0(t)=\frac{t^2}{J^2}.\label{Delta_0-triangle}\eeq 
 The same elimination method provides
$$X_1+X_2=\frac{2J(1+J^2)}{(1-J)^2}.$$
Consequently, the last expression of $\Delta_+(x;t)$ given
in~\Ref{delta-triangle}  reads
\beq \Delta_+(x;t)=1-2x\frac{J(1+J^2)}{(1-J)^2}+x^2J^2,\label{Delta+J}\eeq
and the expression of $S_0(x;t)$ follows, using
Theorem~\ref{complete-small}. Finally, by Theorem~\ref{complete-small}
and~\Ref{delta-triangle},
$$S(x,y;t) =  \frac{\sqrt{\Delta_0(t) \Delta_+(\bx;t)}}
{K(x,y;t)}$$ 
and the announced expression follows, using~\Ref{Delta_0-triangle}
and~\Ref{Delta+J}.

We now apply Theorem~\ref{Sij-small}. The extension of
$\qs(t)$ generated by $\sqrt{\Delta_0(t)}$ and the coefficients of
$\Delta_+(x;t)$ and $\Delta_-(\bx;t)$ is, by~\Ref{Delta_0-triangle},~\Ref{Delta+J} and~\Ref{Jeq}, the field $\qs(J)$.  The polynomials
$A_1(x)=1+x$ and $A_{-1}(x)=1+\bx $ have a   root at $x=-1$, and
$$\sqrt{\Delta_+(-1;t)} =\frac{1+J^2}{1-J},$$
so that the series $S_{i,j}(t)$  belongs to
$\qs(J)$. 
\cqfd
\section{Transcendental examples}
\label{section-transcendental}
\subsection{Extensions of the $(\pm 1, \pm 1)$ model}
We prove in this section that having small {\em down\/} steps is not
enough to ensure algebraicity. We choose the set of steps $\frak S=
\{(1,k), (-1,k), (1,-1), (-1, -1)\}$, with $k\ge 1$. The smallest
positive abscissa of  bilateral walks is 
$$p=\left\{ 
\begin{array}{ll}
1 & \hbox{if } k \hbox{ is even} \\
2 & \hbox{if } k \hbox{ is odd}.
\end{array}
\right.
$$
Moreover, a bilateral walk with $m$ up steps must have $mk$ down
steps, so that its total length is $n=m(k+1)$.
The case $k=1$
was solved explicitly in Theorem~\ref{dim2-diag}, and led to algebraic 
\gfs .  This is no longer true when $k\ge 2$.

\begin{Propo}\label{exemple-transcendant}
 Let $\frak S=
\{(1,k), (-1,k), (1,-1), (-1, -1)\}$, with $k\ge 2$. Let $n$ be a
multiple of $k+1$, say $n=m(k+1)$. If $k$ is even,  assume,
moreover, that $n$ is odd. Then the number of $n$-step walks on the
slit plane ending at $(p,0)$ is
\beq a_{p,0}(n)= \frac{k+1}{m} {n \choose \frac{n+p}{2}}
\sum_{i=0}^{m-1} {{n-1} \choose i} k^{m-i-1}.\label{a-expr}\eeq
The \gf \ $S_{p,0}(t)$ for these walks is D-finite, but
transcendental. The \gfs \  
$L(t), S_0(1;t)$ and $S(1,1;t)$ are not even D-finite.
\end{Propo}
Once again, the proof of this result requires to determine the \gf \
$B(x;t)$ for bilateral walks. We notice that $\frak S$ is a product:
$\frak S= \{-1,1\}\times \{-1,k\}$. This implies that
\beq B(x;t) = \tilde B ((x+\bx)t),\label{Btilde}\eeq
where $\tilde B(t)$ is the \gf \ for bilateral walks with steps in
$\{ (0,-1),(0,k)\}$. 
\begin{Lemma}
Let $U\equiv U(t)$ be the unique power series in $t$ satisfying
$$U=t(1+U)^{k+1}.$$
Then the \gf \ for bilateral walks with steps in
$\{ (0,-1),(0,k)\}$ is
$$\tilde B(t)= \frac{1+U(t^{k+1})}{1-kU(t^{k+1})}.$$
\end{Lemma}
{\bf Proof.}
We begin with the \gf \ $\tilde M(t)$ for bilateral walks that never
go below the $x$-axis. Such walks start with a $(0,k)$ step. By
looking at the first time they visit a point at ordinate $k-1, k-2,
\ldots , 1,0$, we obtain, for their non-commutative    \gf ,
$$\tilde \GM= \epsilon + n ( \tilde \GM s)^k \tilde \GM$$
(the letter $n$ stands for a north step  $(0,k)$, and the letter $s$
for a south step $(0,-1)$). Hence the associated length \gf \ satisfies
\beq \tilde M(t)= 1+t^{k+1} \tilde M(t) ^{k+1}\label{Mtilde eq}\eeq
or, equivalently,
\beq \tilde M(t)= 1+U(t^{k+1}).\label{M-U}\eeq

\begin{figure}[ht]
\begin{center}
\input{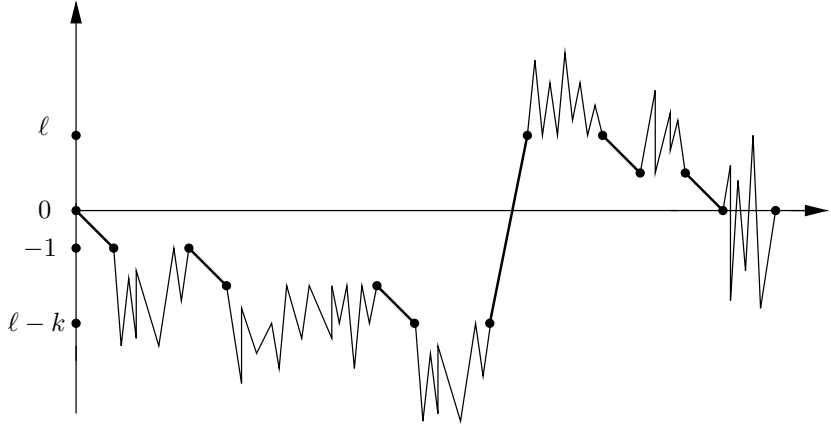}
\end{center}
\caption{Factorisation of a bilateral walk.}
\label{superdyck}
\end{figure}

Now, take a bilateral walk $w$, and look at the first (positive) time $t$
where its ordinate is non-negative. Let $\ell$ denote this ordinate:
then $0\le \ell \le k$ (Figure~\ref{superdyck}). By looking at the
last time before $t$ 
when the ordinate of $w$ is $0, -1,-2,     \ldots , \ell-k$, then at the
first time (after $t$) when the ordinate is $ \ell-1, \ldots ,
1,0$, we obtain:
$$\tilde \GB= \epsilon +\sum_{\ell=0}^k (s\Phi(\tilde \GM))^{k-\ell} n (
\tilde \GM s)^\ell \tilde \GB,$$
where $\Phi(w)$ is the word obtained by reading $w$ from right to
left.
Hence
$$\tilde B(t)= 1+(k+1)t^{k+1}\tilde M(t)^k \tilde B(t),$$
and the announced expression follows from~\Ref{Mtilde eq} and~\Ref{M-U}.
\cqfd

\medskip
\noindent{\bf Proof of Proposition~\ref{exemple-transcendant}.} By
Proposition~\ref{propo-cycle-1}, the 
number of $n$-step walks on the slit plane ending at $(p,0)$ is
\beq a_{p,0}(n) =[x^{p}t^n] \log B(x;t).\label{instance}\eeq 
By~\Ref{Btilde} and the above lemma, we are led to compute, using the Lagrange
inversion formula,
\begin{eqnarray*}
[t^m]\log \frac{1+U}{1-kU} &=& \frac{k+1}{m} [t^{m-1}]
\frac{(1+t)^{m(k+1)-1}}{1-kt} \\
&=& \frac{k+1}{m}\sum_{i=0}^{m-1}  {{(k+1)m-1} \choose i} k^{m-i-1}.
\end{eqnarray*}
The announced expression of $a_{p,0}(n)$ follows.

The series $S_{p,0}(t)$ is D-finite by
Proposition~\ref{propo-cycle-1}. To prove that it is not algebraic, we
are going to study the asymptotic behaviour of its coefficients.
We shall obtain
\beq 
a_{p,0}(n)= \frac{k+1}{\sqrt {2\pi}} \left(
\frac{2(k+1)}{k^{k/(k+1)}}\right)^n  
n^{-3/2} \left( 1- \frac{\sqrt 2 (k-1)}{3\sqrt{\pi k}}\ n^{-1/2}+O(n^{-1})\right), \label{asympt-result}\eeq
and the term in $n^{-2}$ in this expansion  makes it incompatible
with an algebraic \gf \ (see~\cite{Flajolet87}). We start from the
expression of 
$a_{p,0}(n)$ we have just obtained, Eq.~\Ref{a-expr}. An estimate of
the (almost) 
central binomial coefficient can be obtained by standard tools, like
Stirling formula:
\beq \label{stirling} {n\choose (n+p)/2} = \sqrt{\frac 2  \pi} \ 2^n  n^{-1/2} \left(1 +
O(n^{-1})\right) .\eeq
 The rest of the expression of $a_{p,0}(n)$ is the coefficient of $t^m $ in
$\log ({1+U})/({1-kU})$. This series fits perfectly in the framework
of algebraic-logarithmic series developed in \cite{jungen} or
\cite{flajolet-odlyzko}. Recall the equation satisfied by $U$:
\beq P(t,U):=U-t(1+U)^{k+1} =0.\label{eq-U}\eeq
The singularities of $U$ are to be found at the points $t$ such that
$$\frac{\partial P}{\partial u}(t,U)=1-t(k+1)(1+U)^k=0.$$
This implies that $U$ has a unique singularity,  at 
$$t_c= \frac{k^k}{(k+1)^{k+1}},$$
and that $U(t_c)=1/k$. The series $U(t)$ admits an analytic continuation in the domain
$D=\cs\setminus[t_c, +\infty [$. The same holds for $V=
({1+U})/({1-kU})$, because $U=1/k$ forces $t=t_c$.  The logarithm
being an analytic function in $\cs\setminus ]-\infty, 0]$, we have
to mind for the values of $t$ for which $V(t)$ is real and
negative. If $V(t)$ is real, then so are $U(t)$, and $t$ itself
(by~\Ref{eq-U}).
Hence all singularities of $\log V$ lie on the real line. Now $V$
never vanishes on the domain $D$: this would mean that $1+U=0$, so that,
by~\Ref{eq-U}, $U=0$, a contradiction. On $]-t_c,t_c[$, $V(t)$ is real,
given by its series expansion. Moreover, $V(0)=1$, so that by
continuity, $V(t)$ is real positive on all the interval
$]-t_c,t_c[$.
By continuity at $-t_c$, $V(-t_c)$ is real and non-negative. But
$V$ does not vanish, hence  $V(-t_c)>0$. By continuity again, there
exists $t_0>t_c$ such that for $t>-t_0$, $V(t) \not \in \rs ^-$. 

We have thus proved that $\log V$ is not singular in $\cs
\setminus(]-\infty, -t_0] 
\cup [t_c, +\infty[)$. In particular, we can apply to this series the
analysis of singularities of~\cite{flajolet-odlyzko}. From the local expansion
of $U$ around its singularity,
$$U(t)= \frac{1} k -\frac{\sqrt{2(k+1)}}{k\sqrt k}
\sqrt{1-t/t_c}+\frac{2(k+2)}{3k^2} \ (1-t/t_c)+ O( (1-t/t_c)^{3/2}),$$
we obtain
$$\log V(t)=\log  \frac{1+U}{1-kU}=\frac 1 2 \log \frac{k+1}{2k}
+ \frac 1 2 \log \frac 1{1-t/t_c} +\frac{\sqrt 2
(k-1)}{3\sqrt{k(k+1)}} \sqrt{1-t/t_c} +  O( 1-t/t_c).$$
The analytic properties of $\log V$, derived above, allow us to
extract term by term from this estimation the asymptotic behaviour of
the coefficient of $t^m$ in $\log V$:
$$[t^m] \log  \frac{1+U}{1-kU}= \frac 1 {2m} \ t_c^{-m}
-\frac{k-1}{3\sqrt{2\pi k(k+1)}} \ t_c^{-m} m^{-3/2}+ O( t_c^{-m} m^{-5/2}).$$
Eq.~\Ref{asympt-result} is obtained by combination of this estimate
with~\Ref{stirling}.
 The transcendence of $S_{p,0}(t)$ follows.

\medskip

We now turn our attention to the series $L(t)$, $S_0(1;t)$ and
$S(1,1;t)$. By Theorem~\ref{theoreme-general} and
Proposition~\ref{nature-fact}, $L(t)$ is 
D-finite if and only if the coefficient of $x^0$ in $B'/B$ is
algebraic. But
$$[t^{n-1}][x^0]\frac{B'(x;t)}{B(x;t)} = n [t^{n}][x^0] \log B(x;t),$$
and an asymptotic analysis very similar to the one above (one simply
has to replace $p$ by $0$) shows that this coefficient as a term in
$1/n$ in its asymptotic expansion, which prevents the series $[x^0]
B'/B$ from being algebraic. Hence $L(t)$ is not D-finite.

\medskip
The second equation of Proposition~\ref{greg} shows that $S(1,1;t)$ is
D-finite if and only if $S_0(1;t)$ is also D-finite. If  $S_0(1;t)$
was D-finite, then, by the argument of Proposition~\ref{nature-fact},
the series  $[x^>]B'/B$, evaluated at $x=1$, would be algebraic. By
symmetry of $B$ in $x$ and $\bx$, the  series  $[x^<]B'/B$, evaluated
at $x=1$, would also be algebraic. Since $B'/B$ is algebraic, $[x^0]
B'/B$ would finally be algebraic too, so that $L(t)$ would be
D-finite, which we have proved not to be the case.
\cqfd

\subsection{The half-plane with a forbidden half-line}
We conclude this section with a simple model that does not fit perfectly in
the framework studied in this paper, but which we find interesting.

\begin{figure}[htb]
\begin{center}
\epsfig{file=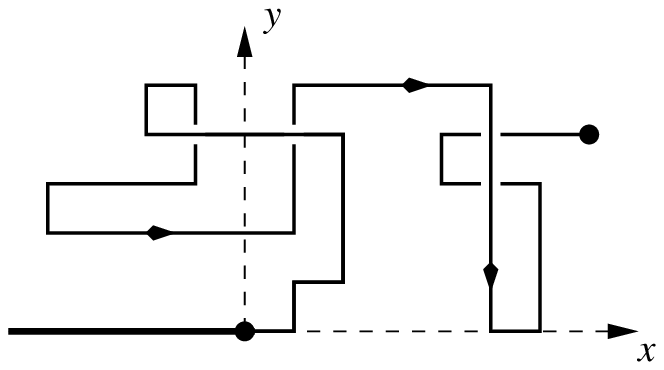}
\end{center}
\caption{A walk on the upper half of the slit plane.}
\label{slit-up}
\end{figure}

Let us go back to the first model studied in Section~\ref{alg-examples}: walks
on the ordinary square lattice that avoid the horizontal half-line ${\cal
H}=\{(k,0), k\le 0\}$.  In 
addition to this constraint, we now force the walks not to visit a
point with a negative ordinate (Figure~\ref{slit-up}). Let $S(x,y;t)$ be
the \gf \ for these walks. Let $S_0(x;t)$ enumerate those that end on
the $x$-axis, and let $L(t)$ be the \gf \ for loops, constrained in
the same manner. We obtain the following counterpart  of Proposition
\ref{greg}. 
\begin{Propo}\label{greg-half}
Let $M(x;t)$ be the \gf \ for bicolored Motzkin walks, given by Lemma~{\em
\ref{lemmeMB}}: 
$$ M(x;t)=\frac{1-t(x+\bx)
-\sqrt{(1-t(x+\bx+2))(1-t(x+\bx-2))}}{2t^2}.$$
The \gfs \ $S_0(x;t)$ and $L(t)$ are related by the following identity:
$$S_0(\bx;t)L(t) S_0(x;t) = M(x;t).
$$
Moreover, the series $S(x,y;t)$ is
related to $S_0(x;t)$ by
$$S(x,y;t)=\frac{S_0(x;t)}{1-tyM(x;t)}.$$ 
\end{Propo}
{\bf Proof.} The proof is a simple adaptation of the argument proving
Proposition~\ref{greg}. We just impose the non-negativity
condition to all walks under consideration. We thus have to replace
ordinary walks (\gf : $1/K(x,y;t)$) by walks that always stay on or
above the $x$-axis.  Clearly, their non-commutative \gf \ $\GW ^+$ 
satisfies $\GW ^+=\GM + \GM\ n \ \GW ^+$, where $\GM$ 
is the non-commutative \gf \ for Motzkin words, so that
 $W^+(x,y;t)=M(x;t)/(1-tyM(x;t))$. Similarly, bilateral walks
(GF: $B(x;t)$) have to
be replaced by Motzkin walks (GF: $M(x;t)$). 
\cqfd
\begin{Propo}
The length \gf \ $S_{1,0}(t)$ for walks on the upper half-plane that
start from $(0,0)$, end at $(1,0)$ and avoid the forbidden half-line
$\{(k,0) : k \le 0\}$ has the following simple expansion:
$$S_{1,0}(t) = \sum_{n \ge 0} \frac{t^{2n+1}}{2n+1} {{2n+1} \choose
n}^2.$$ 
It is D-finite, but transcendental. The series $L(t)$, $S_0(1;t)$ and $S(1,1;t)$ are not D-finite.
\end{Propo}
{\bf Proof.}
The counterpart of Proposition~\ref{propo-cycle-1} is:
$$S_{1,0}(t) = [x] \log M(x;t).$$
This series can be evaluated using the Lagrange inversion
formula (LIF). Let $N(x;t)=t M(x;t)$. Then
$$N= t(1+xN)(1+\bx N),$$
so that $ \log M(x;t)= \log (1+xN)(1+\bx N)$. Using the LIF, we find
$$\log M(x;t)= \sum_{n \ge 1} \frac{t^n}{n} \sum_{k = 0} ^n {{n}
\choose k}^2 x^{2k-n}.$$
Extracting the coefficient of $x$ yields the announced result for
$S_{1,0}$. The coefficient of $t^{2n+1}$ in the  series $S_{1,0}(t)$
grows like 
$4^{2n}/n^2$, up to a multiplicative constant, and this behaviour is
not compatible with an algebraic  \gf \  (see
\cite{Flajolet87}). The rest of the argument follows the same
principle as the end of the proof of
Proposition~\ref{exemple-transcendant}: $L(t)$ is D-finite if and only
if the coefficient of $x^0$ in $M'/M$ is algebraic. But
$$[t^{n-1}][x^0]\frac{M'(x;t)}{M(x;t)} = n [t^{n}][x^0] \log M(x;t)={n
\choose {n/2}} ^2\sim c \ 4^n n^{-1},$$
and this cannot be the asymptotic behaviour of the coefficients of an
algebraic series. The second equation of Proposition~\ref{greg-half} shows that
$S(1,1;t)$ is D-finite if and only if $S_0(1;t)$ is D-finite. As
above, this assumption forces $L(t)$ to be D-finite as well, which is
not the case.
\cqfd

\section{Further directions} \label{directions}

\subsection{Towards a complete classification of the sets of steps} In
view of the results of Sections~\ref{section-algebraicity}
and~\ref{section-transcendental},  we conjecture that the slit plane
model associated with a set of steps $\frak S$ has an algebraic \gf \
if and only if the $x$-axis cannot be crossed without being
effectively visited. This general statement does not apply to
pathological sets $\frak S$ in which no step goes, for instance, right
(or left, or north, or south): the corresponding models are, in
essence, one-dimensional, and yield algebraic series. Our conjecture
means that Theorem~\ref{complete-small} essentially encapsulates all algebraic
cases (if we assume the g.c.d. of all vertical moves to be
$1$). Section~\ref{section-transcendental} suggests that a possible
approach of this conjecture lies in the asymptotic study of the
coefficients of $S_{p,0}(t) =[x^p]\log B(x;t)$. The case where $\frak
S$ is a product set ${\cal H} \times {\cal V}$ looks promising, since
it decouples the problem into one-dimensional problems: in this case,
$$B(x;t)=\tilde B(H(x)t)$$
where $H(x)$ is a polynomial in $x$ and $\bx$ describing the
horizontal moves and $\tilde B(t)$ is the \gf \ for one-dimensional
bilateral walks with steps in $\cal V$. Hence by
Proposition~\ref{propo-cycle-1}, 
$$[t^n]S_{p,0}= [x^p] H(x)^n [t^n] \log \tilde B(t).$$
As counting one-dimensional bilateral walks consists
in evaluating a coefficient of the form $[y^0]V(y)^n$, where $V(y)$ is a Laurent polynomial in
$y$, this case is likely to be attacked by a classical
saddle-point analysis.

\subsection{Walks on the square lattice avoiding a half-line of
rational slope}
In the myriad of sets $\frak S$ one could decide to study in greater  detail,
some are more appealing than the others; for instance, the 
sets $\frak S_{p,q}=\{(1,-p), (-1,p), (0,q), (0,-q)\}$ which stem from
the enumeration of {\em square lattice\/} walks avoiding a half-line
of rational slope.

\begin{figure}[htb]
\begin{center}
\epsfig{file=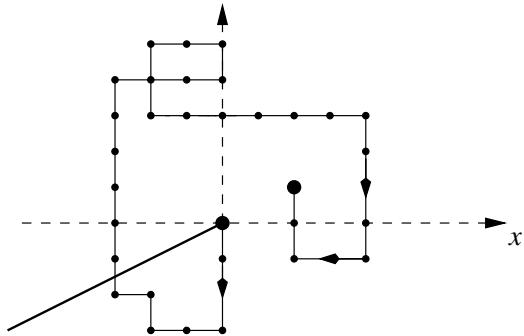}
\end{center}
\caption{This path  crosses the half-line of slope $r=1/2$, but has no
vertex on  it.} 
\label{slope2-fig}
\end{figure}

Indeed, let  $p$ and $q$
be two relatively prime 
integers such that $p\ge 0$ and $q>0$. We consider walks on the
ordinary square lattice that 
start from the origin, but otherwise never meet the integer points of
the half-line ${\cal H}_{p/q}=\{(x,px/q), x\le 0\}$   (Figure~\ref{slope2-fig}). 
 The transformation $(i,j)\mapsto (i,qj-pi)$ shows that counting these
walks is equivalent to counting walks on the slit plane with steps in
$\frak S_{p,q}$. The machinery developed in this paper applies, but
several natural questions could be explored further:
\begin{itemize}
\item
How explicit can the result be made? for slope $2$? for a generic
slope $p/q$?
\item
If the slope of the forbidden half-line is $0$ or $1$, then the \gf \
one obtains is algebraic, as shown by Theorems~\ref{dim2}
and~\ref{dim2-diag}. In accordance with the above conjecture, we are
tempted to believe that these are the only algebraic cases --- but
this needs a proof! A possible approach would be to start from
Proposition~\ref{propo-cycle-1}: the number of $n$-step walks ending
at $(q,p)$ is obtained 
by extracting a coefficient from $\log B(x;t)$, and one could try to
study the asymptotic behaviour of this number.
\item
It is much easier for a square lattice walk to avoid the integers
points of the half-line  ${\cal H}_{100}$ than those of  ${\cal
H}_{0}$. I thought for a while that this might affect significantly
the asymptotic behaviour of the 
number of $n$-step walks avoiding these lines: this
number grows like $4^n n^{-1/4}$ when  ${\cal H}_{0}$ is forbidden,
and I thought that the exponent $-1/4$ could be replaced by a larger
one for  ${\cal H}_{100}$. This intuition was wrong: from the equation
$$S(1,1,t)^2L(t)=\frac 1 {(1-4t)^2B(1;t)}$$
derived from Proposition~\ref{greg}, one can prove, as
Lawler did for the slope $0$ case, that {\em for any rational slope $p/q$,
the number of $n$-step walks avoiding the integer points of  ${\cal
H}_{p/q}$ grows like $4^n n^{-1/4}$} (see~\cite[Section~2.4]{lawler}).
This, of course, only holds for rational slopes:  if the slope is irrational,
 the walks only have to avoid
the point $(0,0)$: their number is directly
related to the return time of a random walk and grows like $4^n/
\log(n)$ (see~\cite[p.~167]{spitzer}). 

\end{itemize}

\subsection{Real avoidance of the half-line}
When we extended the slit plane model from  square lattice walks
to walks with steps in an arbitrary set $\frak S$, we only forbade the
walks to have a vertex on the half-line $\cal H$. The reason for this
choice was simply that the methods we had developed for the original
case worked perfectly for this extension of the model. But what if we also
forbid the walks to {\em cross\/} the half-line? For sets of steps
with small height variations, the model is unchanged, and yields an
algebraic \gf . What about other sets of steps? Can one  of
the methods presented here (or both) be adapted to this new
convention? What can be said of the nature of the associated \gfs ?

\subsection{Higher dimensions}
Let $\frak S$ be a finite subset of $\zs ^d$, with $d\ge 2$. We
consider $d$-dimensional walks that start from the origin of the
lattice $\zs^d$ and have their steps in $\frak S$. The \gf \ for these
walks, counted by their length (variable $t$) and position of their
endpoint (variables $x_1, \ldots , x_d$) is
$$W(x_1, \ldots , x_d;t) = \frac 1 {K(x_1, \ldots , x_d;t)}$$
with
$$K(x_1, \ldots , x_d;t)= 1-t \sum_{(i_1, \ldots , i_d) \in \frak S}
x_1^{i_1} \cdots x_d^{i_d}.$$
A walk is said to be {\em bilateral\/} if it ends on the $x_1$-axis,
that is, on the (one-dimensional) line $\{ (k, 0, \ldots , 0), k\in
\zs \}$. The \gf \ for bilateral walks is
$$B(x_1;t)= [x_2^0 \cdots x_d^0] W(x_1, \ldots , x_d;t).$$
This series  is the iterated diagonal of a rational function: hence it
is always 
D-finite, but usually transcendental when $d\ge 3$.
We consider  walks that never return to the half-line ${\cal H}=\{
(k, 0, \ldots , 0), k\le 0\}$ once they have left their starting
point.  All the machinery developed in
Sections~\ref{section-lawler} and~\ref{section-complete} 
remains valid. In particular, Proposition~\ref{greg} and
Theorem~\ref{theoreme-general} still 
hold, so that the enumeration of these walks boils down to the
explicit factorization of a {\em D-finite\/} (rather than algebraic)
series, namely $B(x_1;t)$. How explicitly can this be done? 

Let us examine the simplest case: walks on the cubic lattice $\zs^3$,
with steps in $\{ (0,0,1), (0,1,0), (1,0,0), (0,0,-1), (0,-1,0),
(-1,0,0)\}$. The projection of a bilateral walk on the $yz$-plane is a
two-dimensional walk that starts and ends at the origin. The number of
such  walks with $2n$ steps is easily seen to be ${{2n}
\choose n}^2$. Hence the \gf \ for bilateral walks on the cubic
lattice is
$$B(x;t) = \sum_{N \ge 0} t^N \sum_{n=0}^{N/2} {N \choose 2n } {{2n}
\choose n}^2 (x+\bx)^{N-2n} .$$
Using the algorithm described in~\cite[Chap.~6]{a=b}, and implemented
by Zeilberger  in the {\sc Maple} package {\sc Ekhad}, one can derive
from the above expression a linear recurrence relation satisfied by
the coefficients of $t^n$ in $B(x;t)$; this recurrence relation can
then be translated into
a following differential equation of order $2$.
%
%
%
What can be said of the
canonical factors of $B(x;t)$? Can they also be described  by differential
equations? Note that the number of
$n$-step walks avoiding $\cal H$ is known to grow like $6^n/\sqrt{\log n}$
(see Eq.~(2.35) in~\cite{lawler}).


\bigskip
\bigskip
\noindent{\bf Acknowledgements.}  I am very grateful to Greg Lawler, who
mentioned to me the  result of his book that was
the starting point of Section~\ref{section-lawler}.
I also wish to thank Gilles Schaeffer, for several reasons: first, he
taught me, 
via his masterful Ph.D. thesis, the charms of the cycle
lemma~\cite{schaeffer-these}; then,  he 
simplified the proof of Proposition~\ref{greg}; and finally, he made
extremely valuable comments on this paper.

\end{document}